\input epsf.sty
\input graphicx
\input  amssym.tex
\def\sqr#1#2{{\vcenter{\hrule height.#2pt              
     \hbox{\vrule width.#2pt height#1pt\kern#1pt
     \vrule width.#2pt}
     \hrule height.#2pt}}}
\def\square{\mathchoice\sqr{5.5}4\sqr{5.0}4\sqr{4.8}3\sqr{4.8}3}
\def\qed{\hskip4pt plus1fill\ $\square$\par\medbreak}

\def\hmap{ {H  } }

\def\cF{{\cal F}}

\def\cX{{\cal X}}

\def\P{{\bf P}}
\def\C{{\bf C}}


\centerline{\bf  Dynamics  of Rational Surface Automorphisms:}

\centerline{\bf  Rotation Domains}
\medskip
\centerline{Eric Bedford\footnote*{Supported in part by the NSF} and Kyounghee Kim}


\bigskip\noindent{\bf  \S0. Introduction. }   Let ${\cal X}$ denote a compact complex surface, and let $f$ be a (biholomorphic) automorphism of ${\cal X}$.  The regular part of the dynamics of $f$ occurs on the Fatou set  ${\cal F}(f)\subset{\cal X}$, where the forward iterates are equicontinuous.   As in [BS, U], we call a Fatou component $U\subset{\cal F}(f)$ a {\it rotation domain of rank $d$} if $f|_U$ generates a (real torus) ${\bf T}^d$-action on $U$.  In dimension 1, rotation domains correspond to Siegel disks or Herman rings, which have a (circle) ${\bf T}^1$- action.  
Here we consider surface automorphisms with the property that the induced map $f^*$ on $H^2({\cal X})$ has an eigenvalue greater than one.   This is equivalent to the condition that $f$ have positive entropy.  

Let us consider generally the possibilities of Fatou sets for surface automorphisms.  If ${\cal X}$ is a complex 2-torus, then an automorphism with positive entropy is essentially an element of $GL(2,{\bf Z})$.  Positive entropy implies that the eigenvalues are $|\lambda_1|<1<|\lambda_2|$, and in this case  the Fatou set is empty.    A second possibility is given by $K3$ surfaces (or certain quotients of them).  Since there is an invariant volume form, the only possible Fatou components are rotation domains.   McMullen [M1] has shown the existence of non-algebraic $K3$ surfaces with rotation domains of rank 2 (see also [O]).  

By Cantat [C], the only other possibilities for compact surfaces with automorphisms of positive entropy are rational surfaces.  In fact,  by [BK2],  the rational case is the most ``frequent.''  By definition, a rational surface is  birationally (or bimeromorphically) equivalent to ${\bf P}^2$, and by a result of Nagata, we may assume that it is obtained by iterated blowups of ${\bf P}^2$.  Rotation domains of rank 1 and 2 have been shown to occur for rational surface automorphisms (see [M2] and [BK1]).  Other maps in this family of rational surface automorphisms were found to have attracting and/or repelling basins (see [M2] and [BK1]).

In this paper we show that positive entropy automorphisms can have large rotation domains.   To describe this, let $\Sigma_0\subset {\bf P}^2$ be the line at infinity.  We will  construct a complex manifolds $\pi:{\cal X}\to{\bf P}^2$ by performing iterated blowups to level 3 over points $\{p_0,\dots,p_{n-1}\}\subset\Sigma_0$.  We let ${\cal F}^1_s$ denote the fiber obtained by blowing up $p_s$, and at level 2 we denote by ${\cal F}^2_s$ the fiber obtained by blowing up a point $q_s\in{\cal F}^1_s$.  We construct a pair $(\hmap,{\cal X})$ with a rotation domain which corresponds to Figure 1:
\proclaim Theorem A.  There is a rational surface ${\cal X}$ with an automorphism $\hmap$ which has positive entropy, and a rotation domain $U\supset \Sigma_0\cup{\cal F}_0^1\cup\cdots\cup {\cal F}_{n-1}^1$.    $U$ is the union of invariant (Siegel) disks on each of which $\hmap$ acts as an irrational rotation.  



By Theorem 3.3 shows that, in addition,  the two fixed points of $\hmap$ in ${\bf C}^2$ are often centers of rank 2 rotation domains.

\smallskip
\epsfysize=1.5in
\centerline{ \epsfbox{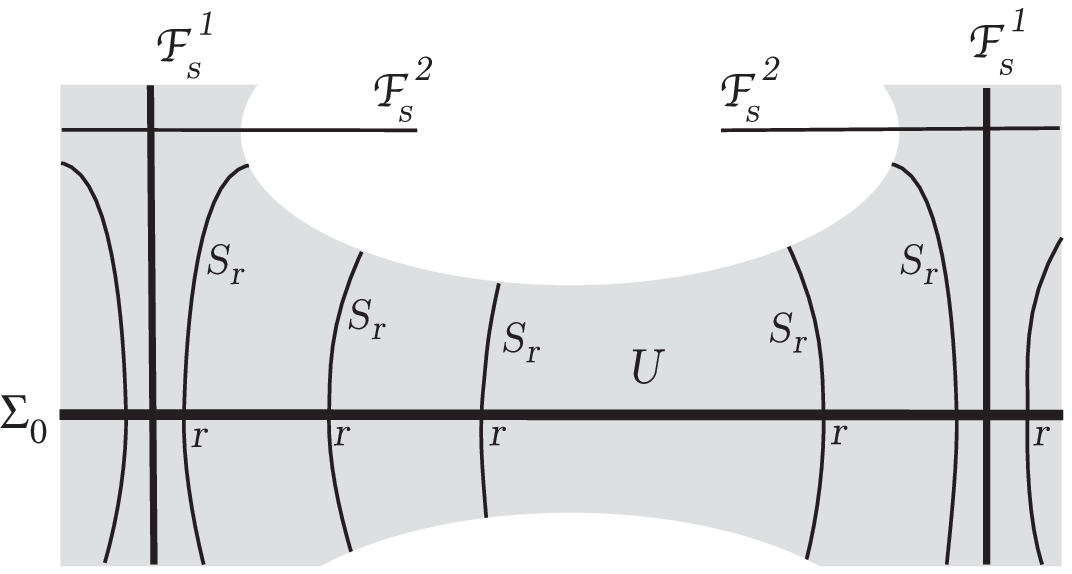}  }

\centerline{Figure 1.  Rotation domain: one parameter family of Siegel disks $S_r$}
\smallskip

For $n,m\ge1$, we define the polynomial
$$\chi_{n,m}(t) = { t(t^{nm}-1)(t^n-2t^{n-1}+1)\over (t^n-1)(t-1)} + 1.\eqno(0.1)$$
If $n\ge4$, $m\ge1$, or if $n=3$,  $m\ge2$,  this is a Salem polynomial, which means that there are real roots $\lambda=\lambda_{n,m}>1>\lambda^{-1}$, and all other roots have modulus one.   We let $\delta$ denote a root of $\chi_{n,m}$ with modulus 1, but which is not a root of unity.  For $1\le j\le n-1$, $(j,n)=1$, we set  $c=2\sqrt\delta\cos(j\pi/n)$, and we define 
$$f(x,y)=(y, -\delta x+cy + y^{-1}).\eqno(0.2)$$
By Theorem 2.3, there is a rational surface $\pi:{\cal X}\to {\bf P}^2$ such that the induced map $f_{\cal X}$ is an automorphism of ${\cal X}$.  We set $\hmap:=f^n$, and will use the pair $(\hmap,{\cal X})$  in Theorems~A and B.

%

Linearization is a useful technique to give the existence of rotation domains, but it is a local technique.  In order to understand the global nature of the Fatou component $U$, we introduce a global model.  We start with the linear map $L=\delta^{n/2}I$ on ${\bf C}^2$, which is scalar times the identity transformation.   $L$ defines a holomorphic map of ${\bf P}^2$ which fixes the line at infinity $\Sigma_0$.  We define a new manifold $\pi:{\cal L}\to{\bf P}^2$ by blowing up ${\bf P}^2$ at 3 levels,  similar to, but different from the procedure used to construct ${\cal X}$.  At each stage, the centers of blowup are fixed points, so $L$ extends to an automorphism of ${\cal L}$, and $(L,{\cal L})$ is our linear model space.  This rotation domain can be linearized on this model space:

\proclaim Theorem B.  There is a domain $\Omega\subset{\cal L}$ and a biholomorphic conjugacy $\Phi:U\to\Omega$ taking $(\hmap,U)$ to $(L,\Omega)$.  In particular, $\hmap$ has no periodic points in $U-\pi^{-1}\Sigma_0$.
 Further, $\pi(\Omega)-\Sigma_0$ is a pseudoconvex, circled domain in ${\bf C}^2$ which is complete at infinity.  

This paper is organized as follows:  \S1 discusses rotation domains generally and global linear models.  \S2 develops a number of the properties of the maps (0.2).  \S3 treats linearization at the fixed points of $\hmap$ in ${\bf C}^2$, which are non-resonant, and $\Sigma_0$ is shown to be in the Fatou set.  In \S4 the resonant fixed points $\cF_s^1\cap\cF_s^2$ are linearized and shown to belong to the same Fatou component as $\Sigma_0$; Theorem A is a consequence of Theorem 4.9.  \S5 gives the global linearization, and Theorem 5.2 yields Theorem B.
\smallskip
\noindent{\it Acknowledgement }  We wish to thank Serge Cantat for several helpful suggestions on this paper.

\medskip

\noindent{\bf \S1. Rotation domains }  
In this section we consider an automorphism $f$ of a general compact, complex manifold ${\cal M}$ of arbitrary dimension.  Recall that the Fatou set consists of all points which have neighborhoods $U'$ such that the restrictions of the forward iterates $\{f^n|_{U'}, n\ge0\}$ form a normal family.   Let $U$ denote an $f$-invariant, connected component of the Fatou set.  We define 
$${\cal G}={\cal G}(U)=\{g:U\to\overline U, g=\lim f^{n_j}\}$$ 
as the set of holomorphic mappings obtained as normal limits of sequences of iterates of $f$.  We say that $U$ is a {\it rotation domain } if ${\cal G}\subset Aut({U})$, i.e.,  every element $g\in{\cal G}$ defines an automorphism of $U$.

\proclaim Proposition 1.1.  If $f$ preserves a smooth volume form, then every Fatou component $U$ is a rotation domain.

\noindent{\it Proof. }  Suppose that $g$ is a normal limit of $f^{n_j}$.  The jacobian of $f^{n_j}$ has modulus one, and so this hold for $g$.  Thus $g$ is an open mapping, so $g(U)\subset U$.  Further, since each $f^{n_j}$ is one-to-one, $g$ also is one-to-one.  Finally, since $g(U)\subset U$, and $g$ preserves volume, we have $g(U)=U$, so $g\in Aut(U)$. \qed

\noindent{\it Remark. } The same argument applies if $f$ preserves a meromorphic volume form, and $U$ has finite volume.

\proclaim Proposition 1.2.  If $U$ is a rotation domain, then ${\cal G}$ is a subgroup of $Aut({\cal M})$.

\noindent{\it Proof. }  If $U$ is a rotation domain, then there is a subsequence $n_j\to\infty$ such that  $\varphi=\lim_{j\to\infty}f^{n_j}\in Aut({\cal M})$.  Passing to a subsequence, we may assume that both $m_j=n_{j+1}-n_j$ and $p_j=n_{j+1}-2n_j$ converge to $+\infty$ as $j\to\infty$.  Passing to further subsequences, we may suppose that there is convergence to limits: $f^{m_j}\to g$ and $f^{p_j}\to h$.  Since $f^{m_j}\circ f^{n_j}=f^{n_{j+1}}$, we see that $g\circ \varphi=\varphi$, so that $g$ is the identity element.  Similarly, $f^{p_j}\circ f^{n_j}=f^{m_j}$, which converges to $h\circ \varphi=g$, which means that $h$ is the inverse of $\varphi$
\qed

The proof of the previous Proposition also allows us to use a characterization of rotation domains similar to one of Fornaess and Sibony [FS]:
\proclaim Proposition 1.3.  A Fatou component $U$ is a rotation domain if and only if there is a subsequence $n_j\to\infty$ such that $f^{n_j}$ converges to the identity uniformly on compact subsets of $U$.

If $U$ is a rotation domain, we have a group action ${\cal G}\times U\to U$.  Since the iterates are a normal family, it follows ${\cal G}$ is a compact group in the compact-open topology. Now we may apply the proof of a Theorem of H. Cartan (as presented, for instance, in Chapter IV of [N]) to conclude:
\proclaim Theorem 1.4.  If $U$ is a rotation domain, then ${\cal G}$ is a compact, abelian, Lie group, and the action of ${\cal G}$ on $U$ is real analytic. 

We let ${\cal G}_0$ denote the connected component of the identity in ${\cal G}$.  Since ${\cal G}$ is a compact, infinite abelian Lie group, ${\cal G}_0$ is a torus of positive dimension $d$.   It is evident that $d\le 2 {\rm\ dim}_{\bf C}{\cal M}$, and we refer to the dimension $d$ as the {\it rank} of the rotation domain $U$. 

\proclaim Theorem 1.5.  If $U$ is a rotation domain, then it is pseudoconvex.

\noindent{\it Proof. }  Pseduconvexity is a local property of the boundary.  The Lie algebra of ${\cal G}$ is generated by holomorphic vector fields.  For a boundary point $p\in\partial U$, we may write a vector field locally in terms of analytic functions $\sum a_j \partial_{z_j}$.  If $U$ is not pseudoconvex, then there will be a coordinate neighborhood on which the $a_j$ have analytic continuations to a larger set.  So the vector field, and thus the torus action,  extends to a larger open set $\tilde U\supset U$.  The larger set $\tilde U$, however, belongs to the Fatou set, which contradicts the fact that $U$ is a Fatou component. \qed

We may regard ${\bf C}^2$ as both an ${\bf R}$-linear and a ${\bf C}$-linear vector space.  Every ${\bf R}$-linear subspace $S\subset{\bf C}^2$ is either complex, or it contains no nonzero complex-linear subspace.  In the second case it is said to be totally real.

\proclaim Theorem 1.6.  Let ${\cal X}$ denote a compact, K\"ahler surface, and let $f$ be an automorphism of ${\cal X}$ with positive entropy.  Then we have $d\le 2$.  If $d=2$, then the generic orbit of ${\cal G}_0$ is a totally real 2-torus, i.e., the tangent space to the orbit is not complex.  If $d=1$, then there is a holomorphic vector field ${\cal V}$, and each orbit of ${\cal V}$ is invariant under $f$.

\noindent {\it Proof. }  Let ${\cal V}_1,\dots, {\cal V}_d$ denote vector fields spanning the Lie algebra of ${\cal G}$.  Let $p$ be a point where their real span has dimension $d$.   We are working inside $T{\cal X}$ which has complex dimension 2, so if $d\ge3$, then we can find a subspace of real dimension 2 in the real span of the ${\cal V}_j$, which is also a 1-dimensional complex subspace of the tangent space $T_p{\cal X}$.   Let  $C$ denote the 2-manifold obtained by moving in the directions spanned by these two vectors.  This will be a 2-torus, which is  a  part of the total ${\cal G}_0$-orbit.  Since the elements of ${\cal G}_0$ are biholomorphic automorphisms, the tangent space to $C$ at each point is a complex submanifold of ${\cal X}$.  We may repeat this argument at any point $p'$ near $p$ and obtain a complex curve $C'$ passing through $p'$.  Since $C$ and $C'$ will be disjoint, we see that $C\cdot C=0$.

On the other hand, since $f$ has positive entropy, there are a cohomology class $\theta_+\in H^{1,1}$ and a $\lambda>1$ such that $f^*\theta_+=\lambda\theta_+$.  By [DF] we must have $\theta_+\cdot\theta_+=0$.   Let $\omega_+$ be a smooth (1,1)-form representing this cohomology class.  We may take the limit $T_+=\lim_{n\to\infty}\lambda^{-n}f^{*n}\omega_+$ and obtain a current which represents the cohomology class $\theta_+$.  However, since the $\{f^n,n\ge0\}$ are a normal family on $U$, it follows that $T_+=0$ on $U$.  We conclude that $\theta_+\cdot C=0$.  However, this makes a 2-dimensional linear subspace of $\{v:v\cdot v=0\}$, which is a contradiction.  Thus $d\le 2$.

Now suppose that $d=2$, and let ${\cal V}_1$ and ${\cal V}_2$ be vector fields which generate the Lie algebra of ${\cal G}_0$.  For generic $p\in{\cal X}$, the span of these vector fields will have real dimension 2.  If the span of these vector fields at a point $p$ is not a complex 1-dimensional subspace of the tangent space $T{\cal X}$, then the ${\cal G}_0$-orbit will be a totally real 2-torus.  Otherwise, if there is an open set where it is complex, we may repeat the argument above.  (Actually, the argument above shows that an orbit of a complex tangency is isolated.)  Thus the generic ${\cal G}_0$-orbit of a point of $U$ must be a totally real 2-torus.

Finally, if $d=1$, then there is a (holomorphic) vector field ${\cal V}$ which generates the Lie algebra.  That is,  ${\cal V}$ generates  a foliation of $U$ by Riemann surfaces, the real part of of ${\cal V}$ generates the action of ${\cal G}_0$.  In particular, each leaf is invariant under ${\cal G}_0$.
\qed

The following remark  shows that torus actions can be linearized even when there is no fixed point.

\proclaim Proposition 1.7.  Suppose that ${\bf T}^2$ is a torus acting on a domain $U$ by biholomorphic automorphisms, and let $z_0\in U$ be a point such that ${\bf T}^2\ni\theta\mapsto\theta\cdot z_0$ is one-to-one, and the orbit of ${\bf T}^2$ is totally real.  Then there is a ${\bf T}^2$-invariant neighborhood $\Omega$ of the orbit ${\bf T}^2\cdot z_0$ and a linearizing map $\Phi:\Omega\to {\bf C}^2$ such that $\Phi(z_0)=(1,1)$, and taking the ${\bf T}^2$-action to the standard ${\bf T}^2$-action on ${\bf C}^2$.

\noindent{\it Proof. }  If we write $\Gamma$ for the orbit ${\bf T}^2\cdot z_0\subset U$, then the equivariance gives a real analytic diffeomorphism from $\Gamma$ to $T:=\{(z_1,z_2)\in {\bf C}^2: |z_1|=|z_2|=1\}$.  Since $\Gamma$ and $T$ are both totally real, this diffeomorphism extends holomorphically to an open set, and this gives the desired linearization.
\qed

In \S3 and \S4 we will use local linearization to show that certain fixed points belong to the Fatou set.   The converse is easier: linearization is always possible at a fixed point inside the Fatou set.    We will recall the statement of an easily proved result about domains in ${\bf C}^k$ (see  [H]), in which $\Phi$ is defined on all of $U$. 
\proclaim Proposition 1.8.  Suppose that $U\subset{\bf C}^k$ is bounded and invariant under a holomorphic map $h$.  Suppose that $z^0=0$ is a fixed point for $h$; let $A$ denote the differential of $h$ at $z^0$, and suppose $A$ is unitary.  Then $\Phi=\lim_{N\to\infty}{1\over N}\sum_{n=0}^{N-1}A^{-n}h^n$ defines a holomorphic map $\Phi:U\to{\bf C}^k$ such that $\Phi\circ h=A\circ \Phi$.

If $U$ is not contained in ${\bf C}^2$, however, we must define a global model if we want to have a global linearization.

\noindent{\bf  Global Linear Model. } Let us give some examples of (zero entropy) maps  which illustrate some possibilities for global rotation.  We start with the linear map on ${\bf P}^2$ which is given as $M[t:x:y] = [t:\mu_1x:\mu_2y]$.    $M$ has 3 fixed points on ${\bf P}^2$: $[1:0:0]$,  $[0:1:0]$, and  $[0:0:1]$.  The multipliers at $[1:0:0]=(0,0)\in{\bf C}^2$ are $\{\mu_1,\mu_2\}$, the multipliers at $[0:1:0]$ (the point where the $x$-axis intersects the line at infinity) are $\{\mu_1^{-1},\mu_2/\mu_1\}$, and at $[0:0:1]$, where the $y$-axis intersects the line at infinity, the multipliers are $\{\mu_2^{-1},\mu_1/\mu_2\}$.  This is shown on the left hand side of Figure 2.  The fixed points are marked; the axes are $X$ and $Y$, and the multipliers at the fixed points are indicated.

{\it Rank 1. } We suppose now that  $|\mu_1|=|\mu_2|=1$, so $M$ generates a torus action if $\mu_1$ and $\mu_2$ are not both roots of unity.  If $\mu_1^{j_1}\mu_2^{j_2}\ne1$ for all $j_1,j_2\in{\bf Z}$, $(j_1,j_2)\ne(0,0)$, then $\mu_1$ and $\mu_2$ are said to be multiplicatively independent, and in this case $M$ generates a ${\bf T}^2$ action on ${\bf P}^2$.

{\it Rank 2. }In the case of multiplicative dependence, we have a ${\bf T}^1$ action.   We may suppose that $\mu_1=t^p$, $\mu_2=t^q$, where $t$ is not a root of unity, and $(p,q)=1$.  Thus $\mu_1^q\mu_2^{-p}=1$. The standard $(p,q)$-action acts on a point $(x,y)\in {\bf C}^2$ according to  ${\bf T}^1\ni\theta\mapsto (e^{ip\theta}x,e^{iq\theta}y)$.   $M$ preserves the curves $\{x^q=cy^p\}$ for any fixed $c\in{\bf C}$.    We  say that $\{\mu_1,\mu_2\}$ are {\it resonant} if $\mu_1^{k_1}\mu_2^{k_2}=\mu_s$ with $s=1$ or $2$, and $k_1,k_2\ge0$ with  $k_1+k_2\ge2$.  There is a special case where the multipliers are $\{1,t\}$, but otherwise in the resonant case, we have a $(p,q)$-action with $pq<0$, which means that only two of the invariant curves pass through the origin.  In the non-resonant case, all of the invariant curves pass through the origin.  

 If $p>q>0$, then the fixed point $[0:1:0]$  will be non-resonant, while the fixed point $[0:0:1]$ will be resonant. 
If $p=q=1$, then $(0,0)$ is non-resonant, but the whole line at infinity $\Sigma_0$ is fixed, with multipliers $\{1,t^{-1}\}$, so all points of the fixed line are resonant.
\bigskip
\epsfysize=1.7in
\centerline{ \epsfbox{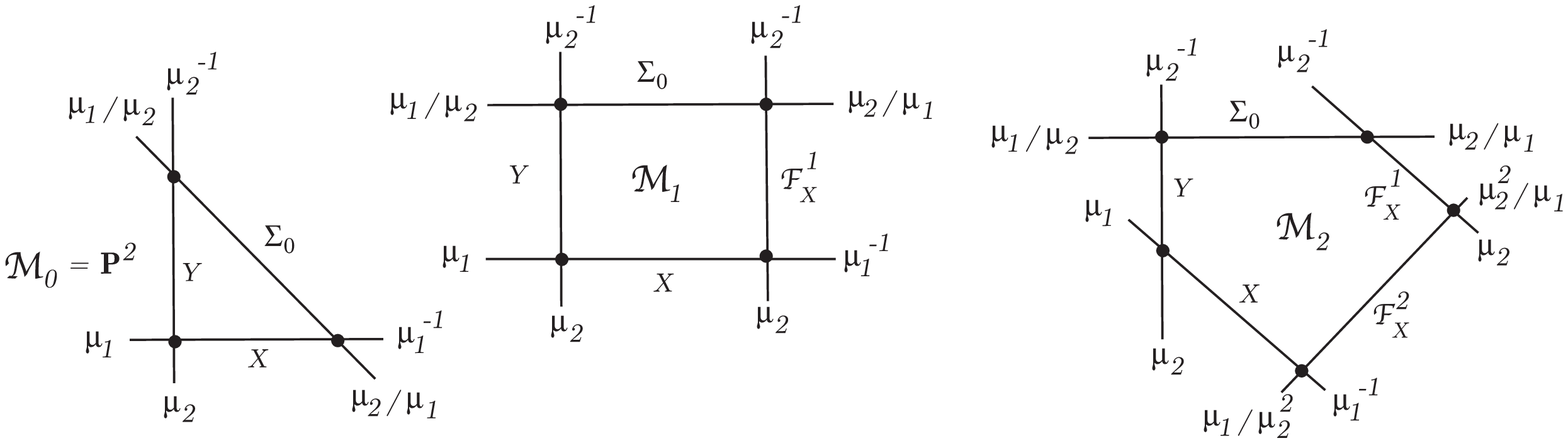}  }

\centerline{Figure 2.  Global linear models ${\cal M}_j$ at stage $j$ of blowup; fixed points and multipliers.  }
\bigskip

{\it Blow up. } Now let $\pi:{\cal Z}\to{\bf P}^2$ be the blow up a fixed point $p$ with multipliers $\{\nu_1,\nu_2\}$, corresponding to directions $X$ and $Y$, respectively.  This will induce the map $\tilde M$ on ${\cal Z}$.  If we denote the resulting exceptional fiber  by $P$, then there will be two fixed points on $P$.  Since $\nu_1$ is the multiplier in the direction $X$, the multipliers at the new fixed point $X\cap P$ will be $\nu_1$ in the direction $X$ and $\nu_2/\nu_1$ in the direction $P$.  Similarly,  the multipliers at  $Y\cap P$ will be $\{\nu_2,\nu_1/\nu_2\}$.   The diagram in the middle of Figure 2 shows the invariant curves and their multipliers after the intersection point $\Sigma_0\cap X$ is blown up; the blowup fibers are is denoted ${\cal F}_X^1$.  The right hand side of Figure 2 shows the space obtained after the further blow up of the point $X\cap {\cal F}_X^1$.   The multipliers at the fixed points are determined by the reasoning described above.   We may repeat this process of blowing up fixed points of $\tilde M$ and obtain a map with an arbitrary number of fixed points.  In the case of rank 1, this produces both resonant and non-resonant fixed points.

\medskip
\noindent{\bf \S2. Rational automorphisms  } Let us imbed ${\bf C}^2$ into ${\bf P}^2$ via the map $(x,y) \mapsto [1:x:y]$. Let $\delta$ be a root of the polynomial $\chi_n$ in (0.1) with $n\ge4$ and $\delta^3\ne-1$, and let $f(x,y)$ be a map of the form (0.2).   In homogeneous coordinates on ${\bf P}^2$, $f$ takes the form
$$f[t:x:y] = [ ty: y^2: -\delta xy  + c y^2 + t^2]. \eqno{(2.1)}$$ 
The exceptional curve for $f$ is $\Sigma_2 = \{ y=0\}$, and $\Sigma_1=\{x=0\}$ is the exceptional curve for $f^{-1}$. 
$$f: \ \  \Sigma_2 \mapsto e_2 = [0:0:1] ,\qquad \qquad f^{-1}: \ \  \Sigma_1 \mapsto e_1 = [0:1:0] $$ Since $f[0:1:w] = [0:1: c- \delta/w]$, the line at infinity $\Sigma_0 =\{t=0\}$ is invariant and $f|_{\Sigma_0}$ is equivalent to the linear fractional transformation $g(w) := c- \delta/w$. Let us set for each $\delta$ such that $\chi_n(\delta) = 0$ , $C_n(\delta) := \{ 2 \sqrt{\delta} \cos ( j\pi/n) :  0 < j < n, (j,n) = 1\}.$
\proclaim Lemma 2.1. If $ c \in C_n(\delta)$ then $f|_{\Sigma_0}$ is periodic with period $n$.

\noindent{\it Proof.} Let $c = 2 \sqrt{\delta} \cos ( j \pi/n)$ for some $j$ relatively prime to $n$. The fixed points of $g$, $w_{\rm fix} = (c \pm \sqrt{c^2- 4 \delta})/2$. It follows that $g'(w_{\rm fix}) = \delta/w_{\rm fix}^2 = e^{2 \pi i j/n}$. \qed

\noindent  It follows that $c \in C_n(\delta)$ if and only if $ g^{n-1} (c) = 0$. Let us use the notation $\omega_s = g^{s-1}(c)$ for $1 \le s \le n-1$, that is $f^s e_2 = [0:1: \omega_s], 1 \le s \le n-1$.

\proclaim Lemma 2.2. Suppose $c \in C_n(\delta)$.  For $1\le j\le n-2$, $\omega_j\omega_{n-1-j}=\delta$.  If $n$ is even, then $\omega_1\cdots \omega_{n-2}=\delta^{(n-2)/2}$.  If $n$ is odd, then we let $\omega_*=\omega_{(n-1)/2}$ denote the midpoint of the orbit.  In this case, we have $\omega_1\cdots \omega_{n-2}=\delta^{(n-3)/2}\omega_*$ and $\omega_*^2=\delta$.

\noindent{\it Proof.} Note that $g^{-1}(w)= \delta/(c-w)$. Since $\omega_{n-1} = 0$, we have $\omega_{n-2} = \delta/c$. It follows that $\omega_1 \omega_{n-2} = c \cdot \delta/c = \delta$. If $\omega_j \omega_{n-1-j} = \delta$ then $\omega_{j+1} = c- \delta/\omega_j$, $\omega_{n-1-(j+(1)} = g^{-1} (\delta/\omega_j) = \delta/(c- \delta/\omega_j)$, and thus $\omega_{j+1} \omega_{n-1-(j+1)} = \delta. $ The Lemma follows by induction on $j$.  \qed

\smallskip
\epsfysize=1.6in
\centerline{ \epsfbox{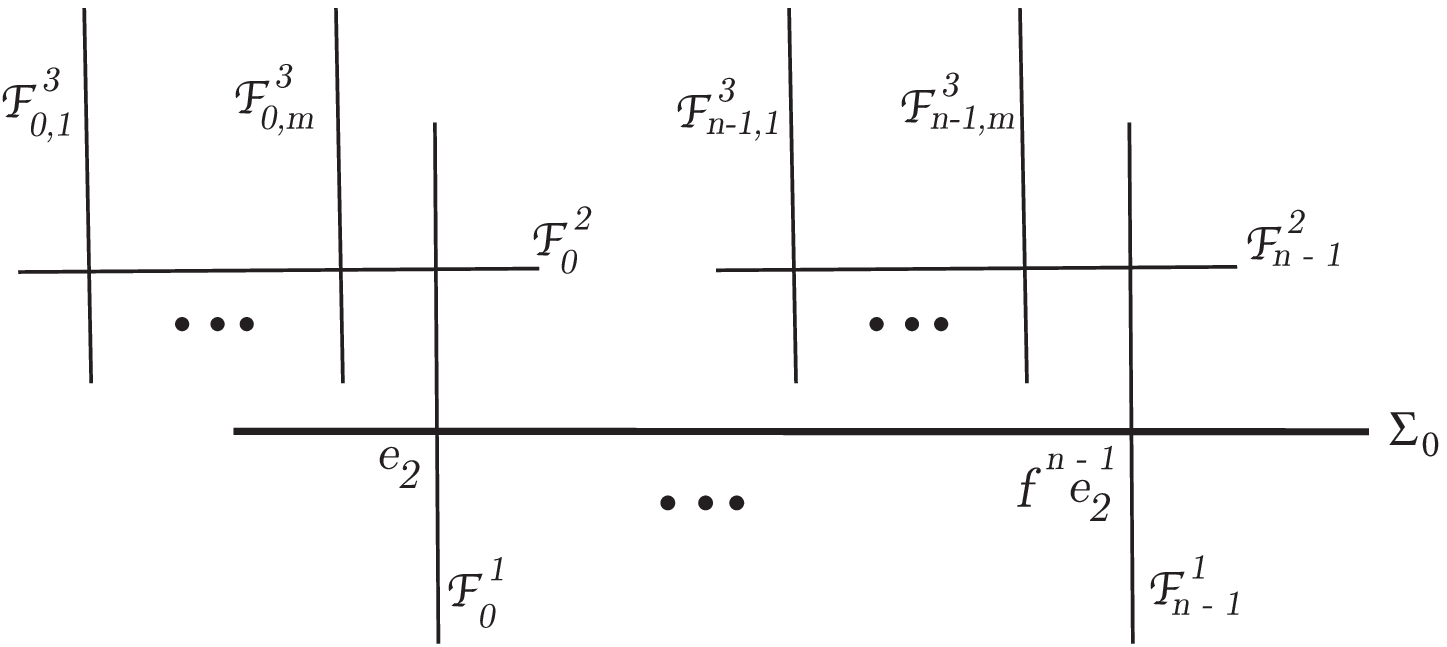}  }

\centerline{Figure 3.  Construction of ${\cal X}$}
\smallskip

For any $\delta$ and any $c \in C_n(\delta)$, we construct the manifold $\pi_1 : {\cal X}^1 \to {\bf P}^2$ by blowing up $n$ points in the line at infinity $f^s e_2, 0 \le s \le n-1$. Let ${\cal F}_s^1:= \pi_1^{-1} (f^s e_2)$ denote the exceptional fibers. For $\cF_0^1$ we will use  $\pi_1(s_1, \eta_1)_0  = [s_1: s_1 \eta_1:1]$ and for $\cF_s^1$, $1 \le s \le n-1$, we use the coordinate chart $\pi_1 (s_1, \eta_1)_s = [s_1: 1: s_1 \eta_1+ \omega_s]$. The induced map $f_{\cX^1}$ maps $\Sigma_2$ to a point in the exceptional fiber $\cF_0^1$ : 
$$f_{\cX^1} [1:x:y] = (s_1, \eta_1)_0= \left ( { y \over 1- \delta x y+ c y^2}, y\right)_0$$
Letting $y \to 0$, we see that $$f_{\cX^1} ( \Sigma_2) = (0,0)_0 = \cF_0^1 \cap \{x=0\}. $$
Similarly we see that $$f_{\cX^1}^{-1} ( \Sigma_1) = (0,0)_{n-1} = \cF_{n-1}^1 \cap \{y=0\}. $$
If we map forward by $f_{\cX^1}$ from $\cF_s^1$ to $\cF_{s+1}^1$,  we have 
$$ f_{\cX^1}\ \ :\ \ \left\{   \eqalign{\ \  &\cF_0^1 \ni (0,\eta_1)_0 \mapsto (0, - \delta \eta_1)_1 \in \cF_1^1\cr 
& \cF_s^1\ni (0,\eta_1)_s \mapsto (0, \delta \eta_1 / \omega_s)_{s+1} \in \cF_{s+1}^1, \quad\ \  1 \le s \le n-2 \cr
& \cF_{n-1}^1 \ni (0,\eta_1)_{n-1} \mapsto (0,\eta_1)_0 \in \cF_0^1.\cr}\right.\eqno{(2.2)}$$
Thus the orbit of the exceptional line $\Sigma_2$ lands at the point of indeterminacy:
$$ f_{\cX^1}:  \cF_s^1 \cap L_s \ni (0,0)_s\mapsto (0,0)_{s+1} \in\cF_{s+1}^1 \cap L_{s+1} \eqno{(2.3)}$$
where $L_0$ is the line $\{ x=0\}$ and $L_s$ is the line $\{ y = \omega_s x\}$, $1 \le s \le n-1$. 

Next, construct $\pi_2: \cX^2 \to \cX^1$ by blowing up the points $(0,0)_s = \cF_s^1 \cap L_s$, for all $0 \le s \le n-1$.  Denote the new fiber by $\cF_s^2$. For $\cF_s^2$ we use local coordinates $\pi_2(\xi_2, x_2)_s = (\xi_2 x_2,x_2)_s = (s_1, \eta_1)_s$ and we have
 $$ f_{\cX^2}\ :\ \left\{   \eqalign{  &(\xi_2,x_2)_0 \mapsto \left({\xi_2\over \xi_2-\delta}\,,\,x_2(-\delta+ \xi_2)\right)_1\cr 
& (\xi_2,x_2)_s \mapsto \left({\omega_2\xi_2 \over \delta x_2^2 \xi_2 + \omega_s(\xi_2+\delta)}\,,\, { x_2 (\delta x_2^2 \xi_2 + \omega_s(\xi_2+\delta)) \over \omega_s( \omega_s + x_2^2 \xi_2)}\right)_{s+1}, \ \  1 \le s \le n-2 \cr
&  (\xi_2,x_2)_{n-1} \mapsto\left({\xi_2 \over \xi_2-\delta+ c x_2^2 \xi_2} \,,\, x_2\right)_0.\cr}\right.\eqno{(2.4)}$$
Thus the induced map $f_{\cX^2}$ the exceptional fibers to the exceptional fibers :
 $$ f_{\cX^2}\ \ :\ \ \left\{   \eqalign{\ \  &\cF_0^2 \ni (\xi_2,0)_0 \mapsto (\xi_2/(\xi_2-\delta),0)_1 \in \cF_1^2\cr 
& \cF_s^2\ni  (\xi_2,0)_s \mapsto (\xi_2/(\xi_2+\delta), 0)_{s+1} \in \cF_{s+1}^2, \quad\ \  1 \le s \le n-2 \cr
& \cF_{n-1}^2 \ni  (\xi_2,0)_{n-1} \mapsto(\xi_2/(\xi_2-\delta),0)_0 \in \cF_0^2.\cr}\right.\eqno{(2.5)}$$
Near $\Sigma_2$ we have $$f_{\cX^2} [1:x:y] = (\xi_2, x_2)_0 = \left( {1 \over 1-\delta x y+ cy^2},y \right)_0.$$
The inverse map is $f^{-1} (x,y) = ( ( cx-y+ {1\over x})/\delta, x)$, which lifts to $$f^{-1}_{\cX^2} (x,y) = (\xi_2, x_2)_{n-1} = \left({\delta \over 1+ cx-cyx},x\right)_{n-1}. $$ 
Thus $f^{-1}_{\cX^2} (\Sigma_1) = ( \delta,0)_{n-1} \in \cF_{n-1}^2$. So in order to have the exceptional curve $\Sigma_2$ land on the point of indeterminacy after $n$ steps, we must have 
$$(f_{\cX^2} )^n \Sigma_2 = f^{-1}_{\cX^2} \Sigma_1 = (\delta,0)_{n-1} \in \cF^2_{n-1} \eqno{(2.6)}$$
Now we use the first line of $(2.5)$ to see that $f_{\cX^2}^2 \Sigma_2 = f_{\cX^2}(1,0)_0 = (1/(1-\delta),0)_1$. We map this point forward by iterating the second part of $(2.5)$ $n-2$ times. The equation $(2.6)$ becomes (projectively): 
$$ \pmatrix{ 1& 0 \cr 1 &\delta}^{n-2} \pmatrix{ 1 \cr 1-\delta } = \pmatrix{\delta \cr 1} $$
From this we see that $(2.6)$ holds if and only if $\delta$ is a root of $\chi_{n,1}$, as defined in $(0.1)$. 

In case (2.6) does not hold, then $f^n_{{\cal X}^2}(\Sigma_2)$ is not indeterminate, and $f_{{\cal X}^2}$ will map it to ${\cal F}^2_0$, and we may map it through the sequence ${\cal F}^2_0\to\cdots\to {\cal F}^2_{n-1}$ again.  Since $\pmatrix{1&0\cr 1&-\delta}^2\pmatrix{1\cr 0} = \pmatrix{1\cr 1}$, we derive an alternative to (2.6):  the condition that $f_{{\cal X}^2}\Sigma_2$ ends up at the point of indeterminacy after $m$ times through this cycle is given (projectively) by
 $$\left( \left( \matrix{1&0 \cr 1& \delta} \right)^{n-2} \left( \matrix{1&0\cr 1&-\delta}\right)^2 \right)^m \left(\matrix{1 \cr 0 } \right) = \left( \matrix{\delta \cr1}\right)$$
This happens exactly when $\delta$ is a root of $\chi_{n,m}$.

We now make the space $\pi_3: \cX^3 \to \cX^2$ by blowing up at the centers $f_{\cX^2} ^{j+1} \Sigma_2 \in \cF_*^2$ for $0 \le j \le nm-1$, and we denote the blowup fiber by $\cF^3_{s,\ell}$ as in Figure 3. We set $\cX : = \cX^3$. Using the similar computation as above we see that the induced map $f_{\cX}$ maps $\Sigma_2$ to the third exceptional fiber : $$ f_{\cX} : \Sigma_2 \ni [t:x:0] \mapsto \left( {x \delta \over t},0\right)_0 \in \cF_{0,1}^3$$ and the mapping from $\cF_{n-1,m}^3$ to $\Sigma_1$ is a local diffeomorphism.  From our construction we conclude:

\proclaim Theorem  2.3.  Let $n,m,j,\delta$ and $f$ be as above.  Then the induced map $f_{\cX}:\cX\to\cX$ is an automorphism.  The exceptional fibers are mapped according to:
$$\eqalign{& \Sigma_0\to\Sigma_0,\ \ \ \ \ \  \ \cF^j_0\to\cF^j_1\to\cdots\to\cF^j_{n-1}\to\cF^j_0,\ \ \  j=1,2\cr
&\Sigma_2\to\cF^3_{0,1}\to\cdots\to\cF^3_{n-1,1}\to\cF^3_{0,2}\to\cdots\to  \cF^3_{n-1,2}\to \cdots\to\cF^3_{0,m}\to\cdots\to \cF^3_{n-1,m}\to\Sigma_1 \cr}$$

Let $S$ denote the span in $Pic(\cX)$ of $\Sigma_0$ and $\cF^j_{s}$, $j=1,2$, $0\le s\le n-1$. The determinant of the intersection matrix $A$ on $S$ is $(3-n) 3^{n-1}$, so ${\rm det} (A) <0$ for all values of $n \ge 4$. Since the dimension of $S$ is $1+ 2n$, which is odd,  and there is at most one positive eigenvalue, all eigenvalues must be strictly negative. It follows that the intersection form is negative definite on $S$, and $Pic(\cX) = S\oplus T$ where $T:=S^\perp$  the orthogonal complement of $S$. Let $\gamma_{s,\ell}$ denote the projection to $T$ of the class $\cF^3_{s,\ell}\in Pic(\cX)$, and let $\lambda_s$ denote the projection of the strict transform of $L_s$ in $\cX$.  

\proclaim Proposition  2.4.  $\lambda_s=\sum_{\ell}\left(-\gamma_{s,\ell}+\sum_{t\ne s}\gamma_{t,\ell}\right)$.
 Thus we may represent the restriction $f_{\cX*}|_{T}$ as
$$\lambda_{n-1}\to\gamma_{0,1}\to\gamma_{1,1}\to\cdots\to\gamma_{n-1,m}\to\lambda_0 =\sum_{\ell}\left(-\gamma_{0,\ell}+\sum_{s\ne0}\gamma_{s,\ell}\right).$$
   The spectral radius of $f_{\cX *}$ is given by the largest zero of the polynomial $\chi_{n,m}$ in $(0.1)$.

\noindent{\it Proof.} We may assume that $s=0$, that is $L_0 = \Sigma_1$. Since $\Sigma_0 = \Sigma_1 \in Pic(\P^2)$, by pulling back by $\pi_1$, we have $\Sigma_1 + \cF_0^1 = \Sigma_0 + \sum_s \cF^1_s \in Pic(\cX^1)$. From $(2.3)$ we see that the center of the blowup for $\cF_0^2$ is $\cF_0^1 \cap L_0$ and there is no centers of blowup of the second blowup fibers in $\Sigma_0$. Thus we have $\Sigma_1 + \cF_0^1+ 2 \cF_0^2= \Sigma_0 + \sum_s( \cF_s^1+ \cF_s^2) \in Pic(\cX^2).$ Pulling back by $\pi_3$ gives 
$$\Sigma_1 + \cF_0^1+ 2 \cF_0^2+ 2\sum_{\ell}\cF_{0,\ell}^3 = \Sigma_0 + \sum_s( \cF_s^1+ \cF_s^2+\sum_{\ell} \cF_{s,\ell}^3) \in Pic(\cX).$$ 
When we project everything to $T = S^\perp$, we have $\lambda_0 = \sum_{\ell}\left(-\gamma_{0,\ell} + \sum_{s\ne 0} \gamma_{s,\ell}\right)$. By Proposition 2.3  we obtain our representation of the restriction $f_{\cX*}|_{T}$. As in [BK2], the spectral radius is given by the restriction $f_{\cX*}|_{T}$, and a direct computation shows that $(0.1)$ is the characteristic polynomial of the transformation defined by the restriction $f_{\cX*}|_{T}$. \qed   
   
\medskip
\epsfysize=1.7in
\centerline{ \epsfbox{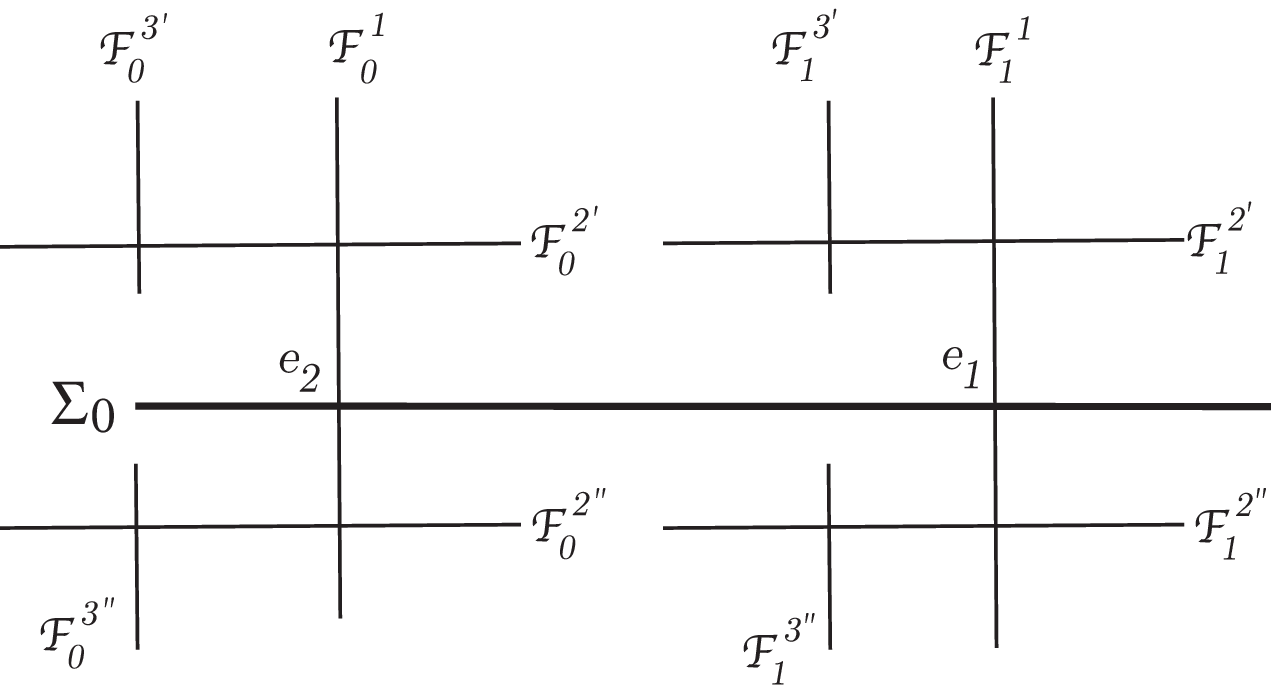}  }

\centerline{Figure 4.  Manifold for Example 2.5.}
\medskip

\noindent{\bf Example  2.5. }  We consider the family of maps given by
$$k(x,y)=\left(y,-x+1+{a\over y}\right).\eqno(2.7)$$
The restriction $k|_{\Sigma_0}$ interchanges $e_1\leftrightarrow e_2$.  As before, we find that $k:\Sigma_2\to e_2$, and $k^{-1}:\Sigma_1\to e_1$, so we blow up the point $e_2$ (resp.\ $e_1$) and denote the resulting fiber as $\cF^1_0$ (resp.\ $\cF^1_1$).  On the new manifold, we have $k:\Sigma_2\to 0\in \cF^1_0$.  The fibers map according to 
$$\cF^1_0\ni\xi\mapsto 1-\xi\in\cF^1_1, \ \ \ \cF^1_1\ni\xi\mapsto\xi\in\cF^1_0.$$
The orbit of $0\in\cF^1_0$ ends up at the point of indeterminacy after going twice around:
$$0\in\cF^1_0\to 1\in \cF^1_1\to 1\in \cF^1_0\to 0\in \cF^1.$$
We blow up this orbit, and label the new fibers so that $\cF^{2'}_0\to\cF^{2'}_1\to\cF^{2''}_0\to\cF^{2''}_1$.  We find that $\Sigma_2$ is still exceptional; it is mapped to $a\in\cF^{2'}_0$.  Finally, we blow up this orbit and obtain an automorphism.

As in the previous case, we consider the invariant subspace $S\subset Pic$, which is generated by $\Sigma_0$ and the blowup fibers up to level 2.  The intersection product restricted to $S$ is negative semidefinite but is not negative definite because it has a zero eigenvalue.  Thus $T=S^\perp$ intersects $S$ in a one-dimensional subspace.  When we compute $k_*|_T$, we find a $3\times3$ Jordan block with eigenvalues of modulus one, so $k_*$ has quadratic growth.

\bigskip\noindent{\bf  \S3. First linearizations }   We will show that $\hmap$ can be linearized at certain fixed points, and thus these fixed points are contained in the Fatou set.  We follow the usual procedure: first we show that there is a formal linearization, which is known to exist  if the multipliers are not resonant.  For the fixed points of $f$ in ${\bf C}^2$, we show that the multipliers are even multiplicatively independent.  The multipliers at the fixed points in $\Sigma_0$ are resonant, but the explicit form of the map there makes the formal linearization clear.   The formal power series is then known to converge if the multipliers satisfy an appropriate ``small divisor'' condition.  In our case, the multipliers are algebraic numbers, and the small divisor property for a pair of multiplicatively independent algebraic numbers is a consequence of the theory of linear forms in logarithms (see Baker [Ba]).  We will also need to linearize the map at certain resonant fixed points, but this is more involved and is postponed to \S4.

The map $f$ has two isolated fixed points in $\C^2$.  Choose one of them and let $\lambda_i, i=1,2$ be the multipliers of the differential.   


\proclaim Lemma 3.1.  The multipliers have modulus $1$ if and only if $ |Re \sqrt{\delta} - 2 \cos (j\pi/n)| \le 1.$

\noindent{\it Proof.}  Direct computation shows that $\lambda_i = - ((1+\delta)/2 - c) \pm \sqrt{ -\delta + ((1+\delta)/2 -c)^2})$, $i=1,2$. Since $|\delta| =1$, we may set $\delta = e^{i\theta}$. Using the expression $c= 2 \sqrt{\delta} \cos(j \pi/n)$ and trig-identities we can see that $$ \lambda_i= \sqrt{\delta} \left[ 2 \cos (j\pi/n) - \cos (\theta/2)  \pm   \sqrt{ ( \cos(\theta/2) - 2 \cos (j \pi/n))^2 -1} \right].$$
Since $|\delta|=1$, two multipliers have modulus $1$ if and only if $( \cos(\theta/2) - 2 \cos (j \pi/n))^2 -1\le 0 $. \qed



\proclaim Lemma 3.2.  The multipliers $\lambda_1$ and $\lambda_2$ are multiplicatively independent. 

\noindent{\it Proof.} Suppose $\lambda_1^{p_1} \lambda_2^{p_2} =1$ for some integers $p_1,p_2 \in {\bf Z}$. Two multipliers $\lambda_1, \lambda_2$ are roots of a polynomial $P(t) = t^2 + (1+ \delta - 2c) t + \delta $. Since $\lambda_1 \lambda_2  = \delta$ and $\delta$ is not a root of unity, we may suppose that $p_1< p_2 \in {\bf Z}$ and $\lambda_1 = \delta ^{p_2/(p_2-p_1)}$ and $\lambda_2 = (1 /\delta )^{p_1 /(p_2-p_1)}$. Let us consider two integers $m,k$ such that $(m,k) =1$ and $p_2/(p_2-p_1) = m/k$ and set $\mu = \delta ^{1/k}$. It follows that $\mu^m$ and $\mu^{k-m}$ are two roots of $P(t)$ and therefore $${1 \over 4} ( \mu^{m-k/2} + \mu^{k/2-m} + \mu^{k/2} +\mu^{-k/2})=  \cos( j\pi/n) \eqno{(3.1)}$$
It is known (see [R, pages 6-7]) that the $n$th Chebyshev Polynomial $T_n$ of the first kind takes on its extrema $\pm 1$ at the points $\cos(j\pi/n)$ for $1\le j \le n-1$.   Thus we have 
$$ T_n(t) +1 = \tau^{(1)}_n(t) \prod_{j:{\rm odd}} ( t- \cos (j \pi/n))\quad {\rm and} \quad T_n(t) -1= \tau^{(2)}_n(t) \prod_{j:{\rm even}} ( t- \cos (j \pi/n))$$
where $\tau^{(i)}_n(t), i=1,2$ are polynomials with no real root. Let us define   $\zeta(t) = ( \mu^{m-k/2} + \mu^{k/2-m} + \mu^{k/2} +\mu^{-k/2})/4$ to be the left hand side of equation $(3.1)$, and let us  set 
$$Q(t) = 4^n t^{kn/2} (T_n( \zeta(t))\pm 1).$$ 
Since $T_n$ is a polynomial of degree $n$, we see that $Q(t)$ is a polynomial with an integer coefficient and $Q(\mu) = 0$. 
Since $\mu^k = \delta$ and $\delta$ is not a root of unity, we conclude that the minimal polynomial of $\mu$ contains exactly one real root, $\mu_*$ outside the unit circle. The minimal polynomial of $\mu$ must divide $Q(t)$ and therefore $Q(\mu_*) =0$ and $\zeta(\mu_*) = \cos (i\pi/n)$ for some $1 \le n$. Since $\mu_*$ is positive real, we have
$$ \zeta(\mu_*) = {1 \over 4} \left( \mu_*^{m-k/2} + {1 \over  \mu_*^{m-k/2}} + \mu_*^{k/2} + {1 \over  \mu_*^{k/2}} \right ) > 1 \ge \cos (j \pi/n) $$ which gives a contradiction.  \qed

\proclaim Theorem 3.3. If $| Re \sqrt{\delta} -2 \cos ( j\pi/n)| \le 1 $, then each of the two fixed points is the center of a rotation domain of rank 2.

\noindent{\it Proof.} By Lemma 3.1 the multipliers $\lambda_1, \lambda_2$ both have modulus $1$ if $| Re \sqrt{\delta} -2 \cos ( j\pi/n)| \le 1 $.  By Lemma 3.2 the two multipliers are multiplicatively independent.  Thus there is a formal power series solution to the linearization equation.  Multiplicative independence means that $T(m_1,m_2):=m_1\log\lambda_1 + m_2\log\lambda_2\ne0$ for all $m_1,m_2\in{\bf Z}$, $(m_1,m_2)\ne(0,0)$.  Since $\cos(j\pi/n)$ and $\delta$ are algebraic, $\lambda_1$ and $\lambda_2$ are algebraic.  It follows from [Ba, Theorem 3.1] that there are $\epsilon>0$ and $\mu<\infty$ such that $|T(m_1,m_2)|\ge\epsilon (|m_1|+|m_2|)^{-\mu}$ for all $m_1,m_2\in{\bf Z}$, $(m_1,m_2)\ne(0,0)$.  This condition is sufficient (see, for instance [P] or [Z]) to show that the formal power series converges in a neighborhood of the origin.  \qed 
\noindent{\it Remark. } We note that for each choice of  $n,m,\delta$, the majority of values of $0<j<n$ satisfy this condition.
\medskip
Now let us turn our attention to the question of linearizing the map $f_{\cX}^n$ in a neighborhood of a point of $\Sigma_0$. To simplify the notation we set $f:= f_{\cX}$ and $\hmap:=f^n$. Let us set 
$$\eqalign { &\lambda:=-\delta^{-n/2}\qquad \quad\qquad\ \  \ {\rm if\ } n:\ {\rm even} \cr & \lambda:= -1/(\delta^{(n-1)/2} \omega_*)\qquad\,{\rm if\ } n:\ {\rm odd}\cr} \eqno{(3.2)}$$  where $\omega_*$ is the midpoint of the orbit defined in Lemma 2.1.

\proclaim Lemma 3.4.  The multipliers of $\hmap$ at $\Sigma_0$ are 1 and $\lambda$.

\noindent{\it Proof.}  Near $\Sigma_0$, $f$ looks like $M:=\pmatrix{0&-\delta\cr 1 & c}$.  We chose $c$ so that the restriction of $\hmap=f^n$ to $\Sigma_0$ will be the identity.  Since $\hmap$ looks like $M^n$ at $\Sigma_0$, we know that $M^n$ should induce the identity map on $\Sigma_0$.  Thus $M^n=\pmatrix{\nu &0\cr 0&\nu}$ is a multiple of the identity matrix, and this means that the multipliers at any point of $\Sigma_0$ will be $\{1,\nu^{-1}\}$.  Since the determinant of $M$ is $\delta$, we conclude that $\nu^2=\delta^n$, or $\nu=\pm \delta^{n/2}$.  It remains to show that the correct sign is the one given in (3.2).

To get the multipliers at $\Sigma_0$ we use  $\tilde\pi_1(\xi_1, t_1)_0  = [t_1 \xi_1:t_1:1]$ for $\cF_0^1$. With this local coordinates, $\{t_1=0\} = \cF_0^1$ and $\{\xi_1=0\} = \Sigma_0$. For $\cF_s^1$, $1 \le s \le n-1$, we use the coordinate chart $\tilde\pi_1 (\xi_1, t_1)_s = [t_1\xi_1: 1: t_1 + \omega_s]$. We also set $h_s =\tilde \pi_1^{-1} \circ f \circ \tilde \pi_1 : ( \xi_1, t_1)_s \mapsto (\xi_1' , t'_1)_{s+1} $ for $ 0 \le s \le n-2$ and $h_{n-1}:  ( \xi_1, t_1)_{n-1} \mapsto (\xi_1' , t'_1)_{0} $. It follows that $f^n = h_{n-1} \circ h_{n-2} \circ \cdots \circ h_0$ in $\Sigma_0 \setminus \{ e_1\}$. Direct computation shows that we have $$ Dh_0 (t_1, 0)_0 = \pmatrix{ -{1/ \delta} & 0 \cr 0 & \star}, \ \  Dh_{n-1} ( t_1, 0) = \pmatrix{ 1 & 0 \cr 0 & \star},$$ and $$Dh_s(t_1,0) = \pmatrix{ { \omega_s /\delta}& 0 \cr 0 & \star} \ 1 \le s \le n-2.$$ It follows that 
$$D\hmap( (t_1,0)_0) = \pmatrix { - ( \omega_1 \cdots \omega_{n-2}) / \delta^{n-1} & 0 \cr 0 & \star} \ \ \ {\rm for\ all\ } (t_1,0)_0 \in \Sigma_0 .$$
Since every point in $\Sigma_0$ is fixed by $\hmap$, we conclude that $\star$, $(2,2)$ entry of the above matrix, is equal to $1$. From Lemma 2.2. we see that $ - ( \omega_1 \cdots \omega_{n-2}) / \delta^{n-1} = \lambda$.  \qed


Let us work in a local coordinate system $(t,\xi)$ near a point $(t=0,\xi=0)\in \Sigma_0=\{t=0\}$.   By Lemma 3.4, the multiplier normal to $\Sigma_0$ at each point is $\lambda$, so we have $\hmap(t,\xi)=(\lambda t + t^2 \star,\xi+t^2 \star)$.   We set $L(t,\xi)=(\lambda t,\xi)$ and consider the equation $\Phi\circ \hmap=L\circ\Phi$, which will give a local conjugacy, conjugating $\hmap$ to $L$.    The form of this particular $\hmap$ is particularly simple (more complicated forms will be considered in the following section), and it is not hard to solve for the higher order terms in the function $\Phi(t,\xi)=(t,\xi) +(t^2\star,t^2\star)$ to obtain a formal solution of this equation.    Since $\lambda$ is algebraic, it satisfies the correct Diophantine condition, and so the series defining $\Phi$ is in fact convergent (see [P], [Ro] or [Ra]):

\proclaim Proposition 3.5. For each $p\in\Sigma_0$, there is a local holomorphic conjugacy $\Phi_p$ at $p$ taking $\hmap$ to a linear map. 

\noindent{\bf  \S4. Linearization at isolated resonant points } Suppose $h$ is a self-map of a general complex $2$-dimensional manifold ${\cal M}$ with a fixed point at the origin. Let $\eta_1$, $\eta_2$ be two resonant multipliers of modulus $1$, that is, $|\eta_i| =1$ and there exist a non-negative integer pair $(a,b) \in {\bf N} \times {\bf N} \setminus \{(0,0)\}$ such that $\eta_1^a \eta_2^b = 1$. It follows that there are infinitely many resonant monomials. Let us define two disjoint sets of monomials spanned by resonant monomials
$$\eqalign{&\breve{{\cal S}_1} =\ {\rm Span} \{x^{j_1}y ^{j_2} :  j_1 = (a/b) j_2 + 1, j_2 \ge 1\} \cr &\breve{{\cal S}_2} =\ {\rm Span}  \{ x^{j_1}y ^{j_2} :  j_2 = (b/a) j_1+1, j_1 \ge 1\}. \cr}$$
Let us also define for each $\{k, \ell \} =\{1,2\}$
$$\eqalign{ & {\cal S}_k =\ {\rm Span} \{x^{j_1} y ^{j_2} : \ (j_1, j_2) \in {\bf N} \times {\bf N} ,   j_k > (a/b) j_\ell + 1 \}, \cr & \hat {\cal S}_k =\ {\rm Span} \{x^{j_1} y ^{j_2} : \ (j_1, j_2) \in {\bf N} \times {\bf N} ,   j_k \ge  (a/b)( j_\ell-1)  \}, .\cr} \eqno{(4.1)}$$

\smallskip
\epsfysize=1.7in
\centerline{ \epsfbox{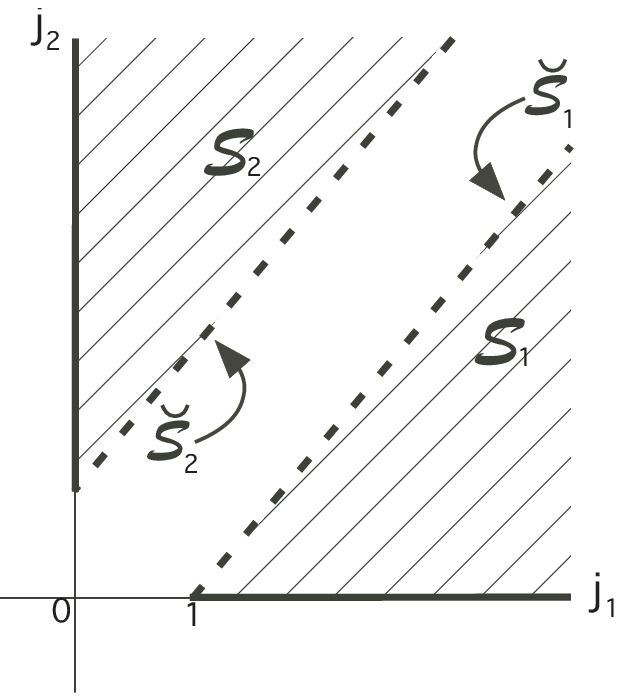}  }

\centerline{Figure 5.  Regions of non-vanishing monomials.}
\smallskip

\proclaim Lemma 4.1. For each $\{k, \ell \} =\{1,2\}$,  ${\cal S}_k$ and $\hat{\cal S}_k$ are closed under multiplication. Furthermore we have 
\itemitem{(a)} For $n\ge 1$, ${\cal S}_k^n = \ {\rm Span} \{x^{j_1} y ^{j_2} : \ (j_1, j_2) \in {\bf N} \times {\bf N} ,   j_k > (a/b) j_\ell + n \} \subset {\cal S}_k$.
\itemitem{(b)}For $n\ge 1$, $\hat{\cal S}_k^n = \ {\rm Span} \{x^{j_1} y ^{j_2} : \ (j_1, j_2) \in {\bf N} \times {\bf N} ,   j_k \ge  (a/b)( j_\ell-n)   \} \subset \hat{\cal S}_k$.
\itemitem{(c)} If $j_1> (a/b) j_2 +1$ then $(x+ {\cal S}_1)^{j_1} (y+ \hat{\cal S}_1)^{j_2} \in {\cal S}_1$
\itemitem{(d)} If $j_1\ge (a/b) (j_2 -1)$ then $(x+ {\cal S}_1)^{j_1} (y+ \hat{\cal S}_1)^{j_2} \in \hat{\cal S}_1$
\itemitem{(e)} If $j_2 > (a/b) j_1 +1$ then $(x+ \hat {\cal S}_2)^{j_1} (y+ {\cal S}_2)^{j_2} \in {\cal S}_2$.
\itemitem{(f)} If $j_2\ge (a/b) (j_1 -1)$ then $(x+ \hat{\cal S}_2)^{j_1} (y+ {\cal S}_2)^{j_2} \in \hat{\cal S}_2$

\noindent{\it Proof.}  Let us suppose $k=1$. Consider an element $s=\sum_{q=1}^m x^{j_{1,q}} y^{j_{2,q}} \in {\cal S}_1$, $m \ge 1$. For $n\ge 1$, $s^n$ is sum of monomials $x^j_1 y^j_2$ where $j_1 = j_{1, q_1}+ \cdots + j_{1,q_n} \ge (a/b) (j_{2, q_1}+ \cdots + j_{2, q_n}) + n = j_2 +n$. It gives the part (a). Similar argument applies for the case (b). In case $j_1 > (a/b) j_2 +1$, using (a) and (b) we have $$(x+ {\cal S}_1)^{j_1} (y+ \hat{\cal S}_1)^{j_2} = (\sum_{i_1} x^{i_1} {\cal S}_1^{j_1-i_1}) \cdot (\sum_{i_2} y^{i_2} \hat {\cal S}_1^{j_2-i_2}) .$$ For each $x^{i_1 + \alpha_1} y^{\alpha_2} \in x^{i_1} {\cal S}_1^{j_1-i_1}$ we have $i_1 + \alpha_1 > (a/b) \alpha_2 + j_1$ and for each $x^{\beta_1} y^{i_2+\beta_2} \in y^{i_2} \hat {\cal S}_1^{j_2-i_2}$, $\beta_1 \ge (a/b) (\beta_2-j_2+i_2)$. It follows that $$ i_1+ \alpha_1 + \beta_1 \ge (a/b) (\alpha_2+ \beta_2 + i_2) + j_1 -(a/b)j_2> (a/b) (\alpha_2+ \beta_2 + i_2) +1.$$ Similarly when  $j_1\ge (a/b) (j_2 -1)$ we have
$$ i_1+ \alpha_1 + \beta_1 \ge (a/b) (\alpha_2+ \beta_2 + i_2) + j_1 -(a/b)j_2> (a/b) (\alpha_2+ \beta_2 + i_2-1).$$ The proof for the case for $k=2$ is essentially identical. \qed

The following Proposition is the direct application of the previous Lemma 4.1.

\proclaim Proposition 4.2. Suppose both $f_i : {\cal M} \to {\cal M}$, $i=1,2$ fix the origin and 
$$ f_i (x, y) \in( \alpha_1^{(i)} x,\  \alpha_2^{(i)} y) +  {\cal S}_1 \times  \hat {\cal S}_1. $$ 
Then we have $$f_1 \circ f_2:(x,y) \mapsto   ( \alpha_1^{(1)}\alpha_1^{(2)}x ,\  \alpha_2^{(1)}\alpha_2^{(2)}y) +  {\cal S}_1 \times  \hat{\cal S}_1. $$ Similarly if $ f_i (x, y) \in ( \alpha_1^{(i)} x,\  \alpha_2^{(i)} y) +  \hat {\cal S}_2\times   {\cal S}_2$ for each $i=1,2$ then so does $f_1 \circ f_2$.

\proclaim Theorem 4.3. If $h$ has a following local expansion: $h(x,y) \in (\eta_1 x, \eta_2 y) +  {\cal S}_1 \times  \hat {\cal S}_1 $ or $h(x,y) \in (\eta_1 x, \eta_2 y) +  \hat{\cal S}_2 \times   {\cal S}_2 $, then there is a formal power series expansion $\Phi$ such that $\Phi\circ h = L \circ \Phi$ where $L(x,y) =(\eta_1 x, \eta_2 y)$. Furthermore if $\eta_1, \eta_2$ are algebraic, then $h$ is linearizable. 

\noindent{\it Proof.} Without loss of generality we may assume that $k=1$. For $j\ge 2$, let us define $\Phi_j: (x,y) \mapsto (x,y) +$ homogeneous polynomials of degree $j$. Since $h$ has no resonant monomials, we can find $\Phi_2$ such that 
$$ \Phi_2^{-1} \circ h \circ \Phi_2 (x,y) =L(x,y) + \tilde h_3$$ 
where order of $\tilde h_3 \ge 3$. From Proposition 4.2. we see that $\tilde h_3 \in {\cal S}_1 \times \hat{\cal S}_1$ and thus it has no resonant monomials. We proceed with an induction on $j$ and find $\Phi_j$ such that $$\Phi_j^{-1} \circ \cdots \circ  \Phi_2^{-1} \circ h \circ \Phi_2 \circ \cdots \circ \Phi_j = L + \tilde h_{j+1}$$ 
where order of $\tilde h_{j+1} \ge {j+1}$ and $\tilde h_{j+1} \in {\cal S}_1 \times \hat{\cal S}_1$. Letting $\Phi = \lim_{n \to \infty} \Phi_2 \circ \Phi_3 \circ \cdots \circ \Phi_n$ we have a power series expansion $\Phi$. 

We define $T(m_1,m_2)=m_1\log\eta_1+m_2\log\eta_2$.  Because of resonances in the multipliers, $T(m_1,m_2)$ can vanish, but by our construction there are no nonvanishing resonant monomials, which means that the coefficient of $x^{m_1}y^{m_2}$ will also vanish when $T(m_1,m_2)$ vanishes.   Since $\eta_1$ and $\eta_2$ are algebraic, by [Ba, Theorem 3.1] we will have $|T(m_1,m_2)|\ge\epsilon(|m_1|+|m_2|)^{-\mu}$ for all values for which $T(m_1,m_2)$, and thus the coefficient of $x^{m_1}y^{m_2}$, does not vanish.  By [P, Z] it follows that the power series of $\Phi$ actually converges, and thus $h$ is linearizable.  \qed

\noindent{\it Remark. }  Another formulation for linearizing  resonant points is given by Raissy  [Ra].
\medskip
\noindent{\bf Example. } Let $\lambda$ is a number of modulus 1 which is not a root of unity and consider the map
$$f(x,y)=(\lambda x,\lambda^{-1}y + x^2+ x^3y).$$
The multipliers at the origin exhibit the resonance corresponding to $\lambda^a(\lambda^{-1})^b=1$ with $a=b=1$.  Thus the monomials in $f$ belong to ${\cal S}_1$, so by Theorem 4.3, $f$
can be formally  linearized at the origin.  If $\lambda$ is algebraic, then the formal conjugacy actually converges and gives a holomorphic linearization of $f$.
\bigskip
We may reformulate Theorem 4.3 to give non-linearizability.
\proclaim Corollary 4.4. Suppose the local expansion is $h= (\eta_1 x, \eta_2 y) + (h^{(1)}, h^{(2)})$. Suppose for $k \in \{1,2\}$, $h_1, h_2 \in \hat{\cal S}_k $ and $h_k \not\in {\cal S}_k$ then $h$ is not linearizable. 

\noindent{\it Proof.} Suppose $n$ is the smallest integer such that $h_k$ has a monomial of order $n$ in $\breve{\cal S}_k$. Using Theorem 4.3. we see that there are $\Phi_2, \dots, \Phi_{n-1}$ such that $$\Phi_{n-1}^{-1} \circ \cdots \circ  \Phi_2^{-1} \circ h \circ \Phi_2 \circ \cdots \circ \Phi_{n-1} = L + \tilde h_{n+1}$$ where $\tilde h_{n}$ has resonant monomials. It follows that there is no formal power series expansion $\Phi$ such that $\Phi\circ h = L \circ \Phi$ and thus $h$ is not linearizable.\qed

Let us apply this discussion to the the map of interest. $\hmap:=f^n$ has isolated fixed points on the exceptional fibers. By Proposition 1.3. we see that exceptional fibers $\cF_i^j$, $0\le i\le n-1$ and $j=1,2$ are all fixed (as a variety) under $\hmap$.  Using local coordinates $(\xi, x)_s$ defined in $(1.6)$,  we rewrite the original map $f$ near $\cF_s^1 \cap \cF_s^2$.  In this coordinate system, $\{\xi=0\} = \cF_s^1$ and $\{x=0\} = \cF_s^2$.

\proclaim Lemma  4.5. The fibers ${\cal F}_s^1\cup{\cal F}_s^2$ are invariant under $\hmap$ for $0\le s\le n-1$.  The multipliers of $H$ at ${\cal F}_s^1\cap{\cal F}_s^2$ are $\lambda^2$ and $1/\lambda$ where $\lambda$ is defined in $(3.2)$.

\noindent{\it Proof. } It suffices to show that the multipliers of $\hmap$  at ${\cal F}_s^1\cap{\cal F}_s^2$ are $\lambda^2$ and $1/\lambda$. Let us rewrite $f$ near ${\cal F}_s^1\cap{\cal F}_s^2$ for $0 \le s\le n-1$ using the local coordinates defined in \S 1. In this coordinate system $(\xi_2, x_2)_s$, we have $\{ \xi_2 = 0\} = {\cal F}_s^1$ and $\{x_2=0\} = {\cal F}_s^2$. Using the expression in $(1.6)$ we see that differential at the origin of each mapping is diagonal and thus we have 
$$ \eqalign{d\hmap|(0,0)_s&= \pmatrix{ -{1/ \delta} & 0 \cr 0 & -\delta} \pmatrix{ {1/ \delta} & 0 \cr 0 & {\delta / \omega_1}} \cdots \pmatrix{ {1 / \delta} & 0 \cr 0 & {\delta / \omega_{n-2}}} \pmatrix{ -{1/ \delta} & 0 \cr 0 &1} \cr &= \pmatrix{ {1 /\delta^n}&0 \cr 0 & - {\delta^{n-1}  /( \omega_1 \cdots \omega_{n-2})}} = \pmatrix{ \lambda^2 &0 \cr 0 & 1/\lambda}.}$$ Last equality in the second line comes from Lemma 2.2.  \qed

Let us recall the mapping restricted on the line at infinity is equivalent to the map $g(w) = c- \delta/w$. Furthermore $c \in C_n(\delta)$ if and only if $g^{n-2}(c) = 0$. 
\proclaim Lemma 4.6. Let $g_n = g^n(c)$ for $n \ge 0$. Then we have $$g_n = c - {\delta \over c} - {\delta^2\over  c^2 g_1} - {\delta^3\over  c^2 g_1^2 g_2} - \cdots - {\delta^n\over  c^2 g_1^2 \cdots g_{n-2}^2 g_{n-1}} .$$

\noindent{\it Proof.}  Note that $c = g_0$. The conclusion is equivalent to
$$g_n-g_{n-1} = - \delta^n/( g_0^2 g_1^2 \cdots g_{n-2}^2 g_{n-1}).$$ 
Since $g_1 = c- \delta/c$, it is easy to see that $g_1 - g_0 = - \delta /g_0$. We proceed by induction on $n$: 
$$g_{n+1} - g_n = g(g_n)- g(g_{n-1})=  - {\delta(g_{n-1}- g_n) \over g_{n-1}g_n}.$$ 
Replacing $g_{n-1}- g_n$ by $  \delta^n/( g_0^2 g_1^2 \cdots g_{n-2}^2 g_{n-1})$ we have the conclusion.  \qed

\proclaim Lemma 4.7. The local expansion of $\hmap$ at the fixed point  ${\cal F}_s^1\cap{\cal F}_s^2$ is given by 
$$\hmap(\xi,x) \in (\lambda^2\xi, {1 \over \lambda} x) + {\cal S}_1 \times \hat{\cal S}_1$$

\noindent{\it Proof.} Using the expression $(1.6)$, we can rewrite the mappings near fixed points  ${\cal F}_s^1\cap{\cal F}_s^2$ as following:
$$ f_{\cX^2}\ :\ \left\{   \eqalign{  &(\xi_2,x_2)_0 \mapsto \left(- {\xi_2 \over \delta}+{\cal S}_1,\,-\delta x_2+ \hat {\cal S}_1\right)_1\cr 
& (\xi_2,x_2)_s \mapsto \left(\sum_{m=0}^\infty {(-1)^m \over \delta \omega_s^m} x_2^{2m} \xi_2^{m+1}  + {\cal S}_1  ,\,    {\delta \over \omega_s} x_2+ \hat {\cal S}_1\right)_{s+1},\ 1 \le s \le n-2\cr
&  (\xi_2,x_2)_{n-1} \mapsto\left(- \sum_{m=0}^\infty { c^m \over \delta^{m+1}} x_2^{2m} \xi_2^{m+1} + {\cal S}_1,\, x_2\right)_0.\cr}\right.$$
Using Lemma 4.1. it suffices to show that the first component of $\hmap$ is in ${\cal S}_1$.  First note that if $j_1= (a/b) j_2 +1$ then for $\xi_2^{\alpha_1}x_2^{\alpha_2} \in{\cal S}_1^{j_1}$ and $\xi_2^{\beta_1}x_2^{\beta_2} \in \hat{\cal S}_1^{j_2}$ we have $\alpha_1 + \beta_1 > (1/2) \alpha_2 +j_1 + (1/2)(\beta_2-j_2)=(1/2) (\alpha_2 + \beta_2) + 1. $ It follows that $${\cal S}_1^{j_1} \hat{\cal S}_1^{j_2} \in {\cal S}_1\quad {\rm for\ } j_1= (a/b) j_2 +1\eqno{(4.2)}$$
Using (4.2) we see that $$f^2 (\xi_2, x_2)_0 =  \left(-\sum_{m=0}^\infty {\delta^{m-2} \over  \omega_1^m} x_2^{2m} \xi_2^{m+1}  + {\cal S}_1  ,\,    -{\delta^2 \over \omega_1} x_2+ \hat {\cal S}_1\right)_{2}$$Using the binomial expansion we can keep track of the coefficient of $x_2^{2n}\xi_2^{n+1}$ in the first coordinate. 
 $$f^3 (\xi_2, x_2)_0 =  \left(-\sum_{m=0}^\infty\delta^{m-3}  \left({\omega_1 \omega_2 +\delta \over  \omega_1^2\omega_2}\right)^m x_2^{2m} \xi_2^{m+1}  + {\cal S}_1  ,\,    -{\delta^3 \over \omega_1\omega_2} x_2+ \hat {\cal S}_1\right)_{3}$$ Again using the binomial expansion we proceed this procedure we see that the first coordinate of $f^n (\xi_2, x_2)_0$ is given by 
$$-\sum_{m=0}^\infty\delta^{m+1-n}  \left(-c + {\delta \over  \omega_1} + {\delta^2\over  \omega_1^2 \omega_2} + \cdots + {\delta^n\over   \omega_1^2  \cdots \omega_{n-3}^2 \omega_{n-2}}\right)^m x_2^{2m} \xi_2^{m+1}  + {\cal S}_1 .$$
Since $c = \omega_1$ and $\omega_j = g^{j-1}(c)$, from Lemma 4.6 the coefficient of $x_2^{2m} \xi_2^{m+1}$ vanishes for all $m\ge 1$ if and only if $  g^{n-2}(c) = 0 $, i.e. $c \in C_n(\delta)$. Since $c$ is chosen such that $ c\in C_n(\delta)$ we have the desired conclusion. \qed

\proclaim Proposition 4.8. There is a holomorphic conjugacy $\Phi$ defined in a neighborhood of ${\cal F}_s^1\cap{\cal F}_s^2$ taking $\hmap$ to a linear map. 

\noindent{\it Proof.} Using Lemma 4.6 and Theorem 4.3, we see that there is a formal expansion of $\Phi$ such that $\Phi\circ \hmap = L \circ \Phi$. Since two multipliers are the power of a root of a Salem polynomial $\chi_n$, two multipliers are algebraic and they are not roots of unity. It follows that $\Phi$ is holomorphic.  \qed

\proclaim Theorem 4.9.  If $|\delta|=1$, there is a Fatou component $U$ which is a rotation domain of rank 1 and which contains  $\Sigma_0\cup{\cal F}^1_0\cup\cdots\cup{\cal F}^1_{n-1}$.  In particular $U$ contains a curve $\Sigma_0$ of fixed points as well as isolated fixed points $\{q_0,\dots,q_{n-1}\}$.

\noindent{\it Proof. }   We will show that there is a neighborhood $U_0$ containing $\Sigma_0\cup\bigcup_s{\cal F}^1_s$ and a unique conjugacy which is tangent to the identity along $\Sigma_0\cup\bigcup_s\{q_s\}$  taking $(\hmap,U_0)$ to $(L,\Phi(U_0))$.   Let $\Phi'$ denote the local conjugacy from Proposition 4.8, which is defined in a neighborhood of $q_s={\cal F}^1_s\cap{\cal F}^2_s$.  Let $\Phi$ denote the local conjugacy from the previous Lemma, which is defined in a neighborhood $U_0$ containing $p_s=\Sigma_0\cap{\cal F}^1_s$.  It suffices to show that these two conjugacies may be analytically continued together to one conjugacy which is defined in a neighborhood of ${\cal F}^1_s$.   

Let us use coordinates $(\xi, x)$ from \S4, such that $q_s=(0,0)$, and ${\cal F}^1_s=\{\xi=0\}$.   The series expressing $\Phi'$ has the form $\sum_{k} \sum_{j\le 2k+1}a_{j,k} x^j \xi^k$, and we may assume that it converges for $\{|x|,|\xi|<1\}$.   Thus if $R<\infty$ and we set  $\epsilon=R^{-2}$, and it follows that the series for  $\Phi'$ converges in  $V'$ which contains $\{|\xi|<\epsilon, |x|<R\}$.

Now let us use coordinates $(s,\eta)$ so that $p_s=(0,0)$, and $\{s=0\}={\cal F}^1_s$.  We may assume that $\Phi$ is defined in $V'':=\{|s|<\epsilon,|\eta|<1\}$.  Choosing $R$ sufficiently large, we may assume that $V'\supset V_0:=\{|s|<\epsilon,{1\over 2}<|\eta|<1\}$.  Now both $\Phi'$ and $\Phi$ conjugate the map $\hmap|_{V_0}$ to the linear map $L(x,y)=(x,\lambda y)$.  It follows that $\tilde \phi(x,y):=\Phi'\circ \Phi^{-1}$ commutes with $L$.  In other words, the second coordinate of $\phi$ satisfies  $\lambda\tilde\phi_2(x,y)=\tilde\phi_2(x,\lambda y)$.  Since $\lambda$ is not a root of unity, we conclude that there is a $c(x)$ so that $\tilde\phi(x,y)=(x,c(x)y)$.  Thus $\tilde \phi$ extends holomorphically to $V''$.  Since we have $\Phi'=\tilde \phi\circ\Phi$, it follows that $\Phi'$ extends analytically to $V'\cup V''$, which is a neighborhood of ${\cal F}^1_s$.

Finally, since $\Phi$ and the extended map $\Phi'$ are both tangent to the identity at $p_s$, they agree in a neighborhood of $p_s$, so they combine to give a conjugacy in a neighborhood of  $\Sigma_0\cup\bigcup_s{\cal F}^1_s$.  \qed

\noindent{\bf \S5.  Global linearization }   In order to give  a global linearization of $\hmap$ on $U$, we define a global model linear model $(L,{\cal L})$ as in \S1.  $L$ is the linear map of ${\bf C}^2$  given by the diagonal matrix $L=diag(\lambda^{-1},\lambda^{-1})$ with $\lambda$ as in (3.2).    The line at infinity $\Sigma_0$ is fixed under $L$, and in the successive blowups,  the multipliers are given as in Figure 2, with $\mu_1=\mu_2=\lambda^{-1}$.  We let  $\Lambda_x$ denote the strict transform of the line $\overline{0x}\subset{\bf P}^2$ in ${\cal M}_j$.  Thus $\Lambda_x$ is invariant.  The manifold ${\cal L}$ is obtained by three stages of blowup; the situation over the point $w_s\in\Sigma_0$ is shown in Figure 6.   The construction of ${\cal L}$ is identical to the first two stages of the construction of ${\cal X}$.   We note that the centers of blowup for the first two stages of blowup in the construction of ${\cal X}$ are in fact fixed points of $L$.  The fixed points of $L|_{{\cal F}_1}$ are $\{p_s,q_s\}$, and the fixed points of $L|_{{\cal F}_2}$ are $\{q_s,r_s\}$.  At the third level we blow up $r_s$, which is fixed by  $L$; to emphasize the difference between ${\cal L}$ and ${\cal X}$, we put a hollow dot in $\cF_s^2$ to denote the $m$ points which were blown up to make $\cX$.   The blowup fiber is denoted by $\tilde{\cal F}^3_s$.  There is a birational map $\iota:{\cal L}\to {\cal X}$ which may be regarded as the identity map in a neighborhood of $\Sigma_0\cup{\cal F}_0^1\cup\cdots\cup{\cal F}^1_{n-1}$.
Making another blowup  in Figure 2, we have:
\proclaim Lemma 5.1.  The local multipliers of $(L,{\cal L})$ at $q_s$ are $\{\lambda^{-1},\lambda^{2}\}$, and at $r_s$ they are $\{\lambda^3, \lambda^{-2}\}$.

\smallskip
\epsfysize=2.2in
\centerline{ \epsfbox{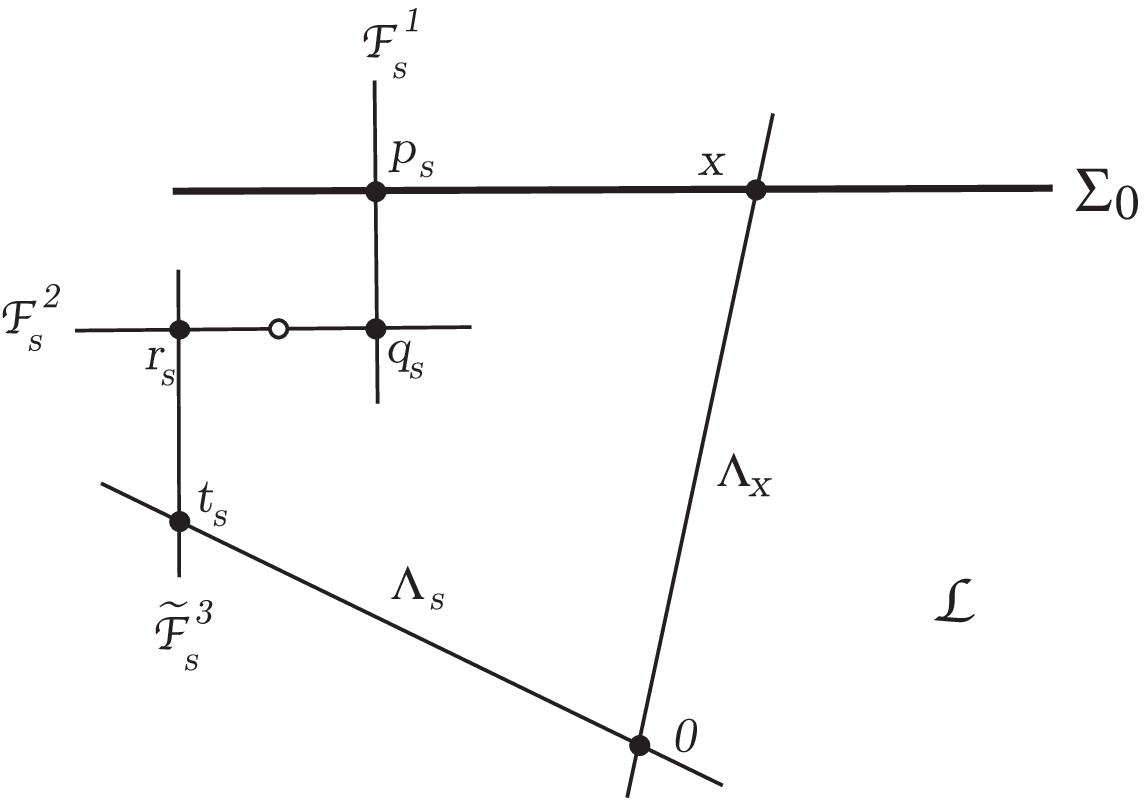}  }

\centerline{Figure 6.  Global linear model ${\cal L}$}
\bigskip

\proclaim Theorem 5.2.  There is a domain $\Omega\subset{\cal L}$ and a holomorphic conjugacy $\Phi:U\to\Omega$ taking $(\hmap,U)$ to $(L,\Omega)$.  In particular, $\hmap$ has no periodic points in $U\cap\pi^{-1}{\bf C}^2$.

One consequence is the following:
\proclaim Corollary 5.3.  The domains $\Omega-\bigcup_s({\cal F}^1_s\cup{\cal F}^2_s\cup\bigcup_{\ell}{\cal F}^3_{s,\ell})$ and  $\Omega':=\pi(\Omega)-\Sigma_0\subset {\bf C}^2$ are pseudoconvex.  Further, $\Omega'$ has the complete-circular property that if $(x,y)\in\Omega'$, and if $\zeta\in{\bf C}$, $|\zeta|\ge1$, then $(\zeta x,\zeta y)\in\Omega'$.

The rest of this section will be devoted to proving Theorem 5.2.

\proclaim Lemma 5.4.  There is a neighborhood $U_0$ containing $\Sigma_0\cup\bigcup_s{\cal F}^1_s$ and a unique conjugacy $\Phi:U_0\to \Phi(U_0)$ which is tangent to the identity along $\Sigma_0\cup\bigcup_s\{q_s\}$  taking $(\hmap,U_0)$ to $(L,\Phi(U_0))$.


\noindent{\it Proof. }  In \S3 we saw that for each point $p'\in\Sigma_0$ there is a local conjugacy $\Phi_{p'}$ between $f$ and the linear map $df_{p'}$ in a neighborhood of $p'$.  This conjugacy was uniquely determined by the condition that its differential is the identity at $p'$.  The construction of $\Phi_{p'}(t,\xi)=(t,\xi)+O_2(t^2)$, moreover, shows that it is tangent to the identity along $\Sigma_0$.  Thus if $\Phi_{p''}$ is another local conjugacy which is defined on an overlapping neighborhood, then $\Phi_{p'}$ and $\Phi_{p''}$ must be continuations of each other.  Now we may use the identification $\iota:{\cal M}\to{\cal X}$ in a neighborhood of $\Sigma_0$ to show that $\Phi$ may be defined in a neighborhood $U_0$ of $\Sigma_0$.  



Let us use coordinates $(\xi, x)$ from \S4, such that $q_s=(0,0)$, and ${\cal F}^1_s=\{\xi=0\}$.   The series expressing $\Phi'$ has the form $\sum_{k} \sum_{j\le 2k+1}a_{j,k} x^j \xi^k$, and we may assume that it converges for $\{|x|,|\xi|<1\}$.   Thus if $R<\infty$ and we set  $\epsilon=R^{-2}$, and it follows that the series for  $\Phi'$ converges in  $V'$ which contains $\{|\xi|<\epsilon, |x|<R\}$.

Now let us use coordinates $(s,\eta)$ so that $p_s=(0,0)$, and $\{s=0\}={\cal F}^1_s$.  We may assume that $\Phi$ is defined in $V'':=\{|s|<\epsilon,|\eta|<1\}$.  Choosing $R$ sufficiently large, we may assume that $V'\supset V_0:=\{|s|<\epsilon,{1\over 2}<|\eta|<1\}$.  Now both $\Phi'$ and $\Phi$ conjugate the map $f|_{V_0}$ to the linear map $L(x,y)=(x,\lambda y)$.  It follows that $\tilde \phi(x,y):=\Phi'\circ \Phi^{-1}$ commutes with $L$.  In other words, the second coordinate of $\phi$ satisfies  $\lambda\tilde\phi_2(x,y)=\tilde\phi_2(x,\lambda y)$.  Since $\lambda$ is not a root of unity, we conclude that there is a $c(x)$ so that $\tilde\phi(x,y)=(x,c(x)y)$.  Thus $\tilde \phi$ extends holomorphically to $V''$.  Since we have $\Phi'=\tilde \phi\circ\Phi$, it follows that $\Phi'$ extends analytically to $V'\cup V''$, which is a neighborhood of ${\cal F}^1_s$.

Finally, since $\Phi$ and the extended map $\Phi'$ are both tangent to the identity at $p_s$, they agree in a neighborhood of $p_s$, so they combine to give a conjugacy in a neighborhood of  $\Sigma_0\cup\bigcup_s{\cal F}^1_s$.  \qed
%


Now let us set $\Sigma_0':=\Sigma_0-\{p_0,\dots,p_{n-1}\}$, and fix $x\in\Sigma_0'$.  The restriction $\Phi^{-1}|_{\Phi(U_0)\cap\Lambda_x}$ is analytic in a neighborhood of $x$, and we let $\omega_x\subset\Lambda_x$ denote a maximal domain such that $\Phi^{-1}$ has an analytic continuation to a map $\psi_x:\omega_x\to U$.  Since $\psi_x$ preserves the circle action, $\omega_x\subset\Lambda_x$ is a disk centered at $x$.

\proclaim Lemma 5.5.  For $x\in\Sigma_0'$, $\omega_x$ is a proper sub-disk of $\Lambda_x-\{0\}$.

\noindent{\it Proof. }  First we observe that $\omega_x$ cannot be all of $\Lambda_x$.  Otherwise, $\psi_x(\omega_x)=\psi_x(\Lambda_x)$ is an algebraic.  However, since $x$ is fixed, this curve is invariant.  But the only invariant curves are $\Sigma_0$, ${\cal F}^1_s$ and ${\cal F}^2_s$, which are not $\psi_x(\omega_x)$.

Now we suppose that $\omega_x=\Lambda_x-\{0\}$.  We have seen that $\psi_x(\omega_x)$ cannot be contained in an algebraic curve.  Now let us define
$$A_{r}:={1\over \#_{r}}\int_0^r[\psi_x(D_t)]\,{dt\over t}$$
where $D_t=\Lambda_x\cap\{|(x,y)|>1/t\}$, and  $\#_{r}$ denotes the area of $\psi_x(D_t)$.  
Passing to a subsequence $r_j\to\infty$, we may construct an Ahlfors current $A$.  Since $x$ is fixed, and $\psi_x$ commutes with the circle action, we have a current satisfying  $\hmap_*(A)=A$.  The corresponding class $\{A\}\in Pic({\cal X})$ is fixed under $\hmap_*$.  Further, a property of the Ahlfors current is that it has nonnegative self-intersection $\{A\}^2\ge0$ (see [Br]).  By \S2, the only elements in $Pic({\cal X})$ are in $S$, and the intersection form is negative definite on $S$.   This contradiction shows that $\omega_x$ must be a proper sub-disk of $\Lambda_x$.  \qed




Recall that by \S1, there is a holomorphic vector field $V$ on $U$ such that $Re(V)$ gives the ${\bf T}^1$ action on $U$.  We let ${\cal V}$ denote the foliation of $U$ which is the complexification of the ${\bf T}^1$ action; that is, ${\cal V}$ consists of the complex leaves of $V$.  For each $x\in\Sigma_0'$, $\psi_x(\omega_x)$ is contained in a leaf of ${\cal V}$, and by  maximality, it is the whole leaf.

\proclaim Lemma 5.6.  For $x\in\Sigma_0'$, the map $\psi_x:\omega_x\to U$ is proper.

\noindent{\it Proof. }   If $K\subset U$ is compact, there exists $\eta>0$ such that for each $y_0\in K$, the leaf of ${\cal V}$ passing through $y_0$ has inner radius at least $\eta$.  We may assume that $\eta$ is less than the distance from $K$ to $\partial U$ and  let $\tilde K$ denote the closure of an $\eta$-neighborhood of $K$.   Thus $\tilde K$ is a compact subset of $U$.  The circle action on $\Lambda_x$ is generated by the vector field $i\zeta{\partial\over\partial\zeta}$, and $\psi_x$ maps this to a constant multiple of $V$.  Since $V$ is bounded on $\tilde K$, it follows that there is a constant $M$ such that the differential of $\psi_x$ at $\zeta_0$ is bounded by $M$ for all $\zeta_0\in\omega_x$ such that $\psi_x(\zeta_0)\in\tilde K$.  It follows that if $\zeta_0\in\psi_x^{-1}K$, then $\psi_x$ extends to the disk of radius $\eta/M$ centered at $\zeta_0$.  Thus the distance of $\psi_x^{-1}K$ to $\partial\omega_x$ is at least $\eta/M$, so $\psi_x$ is proper.  \qed

\noindent{\it Proof of Theorem 5.2.  }  First we assume that  ${\cal F}^2_s\not\subset U$, and define
$$\Omega=\bigcup_{x\in\Sigma_0'}\omega_x \cup\bigcup_{s}({\cal F}_s^1\cup \omega_{q_s}).$$
To see that this is an open set, let us fix $x'\in\Sigma_0'$.  Recall the  associated foliation ${\cal V}$ on $U$, and for $x\in\Sigma_0'$, let $\Gamma_x$ denote the leaf of ${\cal V}$ passing through $x$.  By Lemma 5.6, $\psi_x(\omega_x)=\Gamma_x$, and a property of foliations is that $\Gamma_{x'}\subset \liminf_{x\to x'}\Gamma_x$.  Since $\psi_x$ is equicontinuous, we conclude that $\omega_{x'}\subset\liminf_{x\to x'}\omega_x$, which means that $\Omega$ is open.

We define $\Psi:\Omega\to U$ by setting $\Psi|_{\omega_x}:=\psi_x|_{\omega_x}$.  This map agrees with $\Phi^{-1}$ on a neighborhood of $\Sigma_0\cup{\cal F}^1_0\cup\cdots\cup {\cal F}^1_{n-1}$, and so it is holomorphic there.  Since it is holomorphic when restricted to each $\omega_x$, it follows that $\Psi$ is holomorphic on $\Omega$.  Further, we have seen that $\Psi$ is injective on each $\omega_x$, and disjoint disks $\omega_x$ are mapped to disjoint leaves $\Gamma_x$, so $\Psi$ is injective.

By the semicontinuity of $x\mapsto \omega_x$ and the properness of $\psi_x|_{\omega_x}$, we see that $\Psi$ is proper.  Thus it has a mapping degree, which must be one, and thus $\Psi$ is a biholomorphic conjugacy.

Now if ${\cal F}^2_s\subset U$, then we define
$$\Omega=\bigcup_{x\in\Sigma_0'} \omega_x \cup\bigcup_{s} ({\cal F}^1_s\cup{\cal F}^2_s\cup \omega_{r_s}).$$
All of the previous arguments apply in this case, except that we need to show that $\Psi$ is holomorphic in a neigborhood of $r_s$.  This is similar to the proof of Theorem 4.9.  We may choose coordinates so that ${\cal F}^2_s=\{x=0\}$, and $r_s=(0,0)$.  We may suppose that the map $\Phi$ is holomorphic on the set $\{|x|>1, |y|<1\}$.  Further, since $r_s$ is in the Fatou set, we know that $H$ can be linearized in a neigborhood of $r_s$.  Thus we have another map $\Phi'$ conjugating $H$ to its linear part, which by Lemma 5.1 is $diag(\lambda^{-2},\lambda^3)$.  We may assume that $\Phi'$ is analytic on $\{|x|<2, |y|<1\}$, and that $\Xi:=\Phi'\circ \Phi^{-1}=\sum a_{i,j}x^iy^j$ commutes with this linear map on the set $\{1<|x|<2,|y|<1\}$.  We then have that the first coordinate is $\lambda^{-2}\Xi^{(1)}(x,y) = \Xi^{(1)}(\lambda^{-2}x,\lambda^3y)$, from which we conclude that $\lambda^{-2}=\lambda^{-2i+3j}$ for all nonvanishing coefficients $a_{i,j}$.  Since we have $j\ge0$, and $\lambda$ is not a root of unity,  it follows that we must have $i\ge1$.  Looking at the second coordinate, we get $\lambda^3=\lambda^{-2i+3j}$, so in this case, we cannot have $i<0$.  It follows that all exponents $i,j$ in $\Xi$ are positive, so $\Xi$ is analytic in $\{|x|<2,|y|<1\}$.  Thus we conclude that $\Xi$, and thus $\Psi=\Phi^{-1}$ extends holomorphically through $r_s$.
\qed

\centerline{\bf References}

\item{[Ba]} A. Baker, {\sl Transcendental Number Theory}, Cambridge U. Press, 1990.


\item{[BK1]}  E. Bedford and K. Kim,  Dynamics of rational surface automorphisms: Linear fractional recurrences, J  Geom Anal (2009) 19: 553--583.

\item{[BK2]}  E. Bedford and K. Kim, Continuous Families of rational surface automorphisms with positive entropy,  arXiv:0804.2078

\item{[BS]}  E. Bedford and J. Smillie, Polynomial diffeomorphisms of $C\sp 2$. II. Stable manifolds and recurrence. J. Amer. Math. Soc. 4 (1991), no. 4, 657--679.


\item{[Br]}  M. Brunella, Courbes enti\`eres et feuilletages holomorphes, L'Enseignement mat\'ematique, t. 45 (1999), 195--216.

\item{[C]}  S. Cantat,  Dynamique des automorphismes des surfaces projectives complexes.  C. R. Acad. Sci. Paris S\'er. I Math. 328 (1999), no. 10, 901--906.



\item{[FS]} J.-E. Forn\ae ss and N. Sibony,  Classification of recurrent domains for some holomorphic maps, Math. Ann. 301 (1995), 813--820.

\item{[H]}  M. Herman, Recent results and some open questions on Siegel's linearization theorem of germs of complex analytic diffeomorphisms of $C\sp n$ near a fixed point. VIIIth international congress on mathematical physics (Marseille, 1986), 138--184, World Sci. Publishing, Singapore, 1987.

\item{[M1]}  C. McMullen, Dynamics on $K3$ surfaces: Salem numbers and Siegel disks. J. Reine Angew. Math. 545 (2002), 201--233. 

\item{[M2]}  C. McMullen,  Dynamics on blowups of the projective plane, Publ. Math. Inst. Hautes \'Etudes Sci. No. 105 (2007), 49--89.

\item{[Na]}  R. Narasimhan,  {\sl Several Complex Variables} Reprint of the 1971 original. Chicago Lectures in Mathematics. University of Chicago Press, Chicago, IL, 1995.


\item{[O]}  K. Oguiso,  The third smallest Salem number in automorphisms of $K3$ surfaces.  arXiv:0905.2396 

\item{[P]} J. P\"oschel, On invariant manifolds of complex analytic mappings near fixed points. Exposition.  Math. 4 (1986), No. 2, 97--109.

\item{[Ra]}  J. Raissy, Linearization of holomorphic germs with quasi-Brjuno fixed points,   arXiv:0710.3650

\item{[R]} T.J. Rivlin,  Chebyshev Polynomials, John Wiley \& Sons,  1990

\item{[Ro]}  F. Rong, Linearization of holomorphic germs with quasi-parabolic fixed points. Ergodic Theory Dynam. Systems 28 (2008), no. 3, 979--986. 


\item{[U]}  T. Ueda,   Critical orbits of holomorphic maps on projective spaces. J. Geom. Anal. 8 (1998), no. 2, 319--334.

\item{[Z]} E. Zehnder, A simple proof of a generalization of a theorem by C. L. Siegel. Geometry and topology (Proc. III Latin Amer. School of Math., Inst. Mat. Pura Aplicada CNPq, Rio de Janeiro, 1976), pp. 855--866. Lecture Notes in Math., Vol. 597, Springer, Berlin, 1977.

\bigskip
\rightline{Indiana University}

\rightline{Bloomington, IN 47405}

\rightline{\tt bedford@indiana.edu}

\bigskip
\rightline{Florida State University}

\rightline{Tallahassee, FL 32306}

\rightline{\tt kim@math.fsu.edu}

\bye